\documentclass{amsart}

\usepackage[top=35truemm,bottom=30truemm,left=30truemm,right=30truemm]{geometry}

\usepackage[greek,russian,english]{babel}

\usepackage{amssymb,amsmath,amsthm}    
\usepackage{cite}   
\usepackage{comment}
\usepackage[dvipdfmx]{graphics,color}

\everymath{\displaystyle}

\newcommand\nc{\newcommand}

\theoremstyle{definition}

\newtheorem{theorem}{Theorem}[subsection]
\newtheorem{prop}[theorem]{Proposition}
\newtheorem{importnota}[theorem]{Important Notation}

\newtheorem{notation}[theorem]{Notation}
\newtheorem{defin}[theorem]{Definition}
\newtheorem{caution}[theorem]{Caution}
\newtheorem{remark}[theorem]{Remark}
\newtheorem{lemma}[theorem]{Lemma}
\newtheorem{construction}[theorem]{Construction}
\newtheorem{corollary}[theorem]{Corollary}
\newtheorem{example}[theorem]{Example}
\newtheorem{conclusion}[theorem]{Conclusion}
\newtheorem{triviality}[theorem]{Triviality}
\newtheorem{proto}[theorem]{Prototype Quasifibration}
\newtheorem{cauex}[theorem]{Cautionary Example}
\newtheorem{propositiondef}[theorem]{Proposition-Definition}
\newtheorem{subth}{Nuisance}[theorem]
\newtheorem{conjecture}[theorem]{Conjecture}
\newtheorem{sidest}[theorem]{Side Story}
\newtheorem{miniexample}[theorem]{Example}

\nc\tri{\begin{triviality}}
\nc\side{\begin{sidest}}
\nc\conj{\begin{conjecture}}
\nc\prodef{\begin{propositiondef}}
\nc\prt{\begin{proto}}
\nc\lem{\begin{lemma}}
\nc\sblm{\begin{sublemma}}
\nc\propo{\begin{prop}}
\nc\thm{\begin{theorem}}
\nc\cor{\begin{corollary}}
\nc\defi{\begin{defin}}
\nc\sthm{\begin{subth}}
\nc\exm{\begin{example}}
\nc\miniexm{\begin{miniexample}}
\nc\plm{\begin{problem}}
\nc\rem{\begin{remark}}
\nc\subrmk{\begin{subremark}}
\nc\nota{\begin{notation}}
\nc\cau{\begin{caution}}
\nc\imn{\begin{importnota}}
\nc\cax{\begin{cauex}}
\nc\con{\begin{construction}}
\nc\cnc{\begin{conclusion}}
\newcommand{\pf}[1][]{\begin{proof}[Proof {#1}.]}

\nc\elem{\end{lemma}}
\nc\esblm{\end{sublemma}}
\nc\eside{\end{sidest}}
\nc\econj{\end{conjecture}}
\nc\eprodef{\end{propositiondef}}
\nc\eprt{\end{proto}}
\nc\ethm{\end{theorem}}
\nc\ecor{\end{corollary}}
\nc\edefi{\end{defin}}
\nc\esthm{\end{subth}}
\nc\epropo{\end{prop}}
\nc\etri{\end{triviality}}
\nc\eexm{\end{example}}
\nc\eminiexm{\end{miniexample}}
\nc\erem{\end{remark}}
\nc\subermk{\end{subremark}}
\nc\eplm{\end{prblm}}
\nc\ecau{\end{caution}}
\nc\ecax{\end{cauex}}
\nc\eimn{\end{importnota}}
\nc\enota{\end{notation}}
\nc\econ{\end{construction}}
\nc\ecnc{\end{conclusion}}
\nc\epf{\end{proof}}

\newcommand{\bZ}{\mathbb{Z}}

\newcommand{\bZge}[1]{\mathbb{Z}_{\ge{#1}}}
\newcommand{\bQ}{\mathbb{Q}}
\newcommand{\bR}{\mathbb{R}}
\newcommand{\bC}{\mathbb{C}}
\newcommand{\bP}{\mathbb{P}}
\newcommand{\bI}{\mathbb{I}}

\newcommand{\cH}{\mathcal{H}}
\newcommand{\cA}{\mathcal{A}}
\newcommand{\cU}{\mathcal{U}}
\newcommand{\cV}{\mathcal{V}}
\newcommand{\cD}{\mathcal{D}}

\newcommand{\diff}{d}

\newcommand{\myzeta}[2]{\zeta_{F_2}{\begin{pmatrix}{#1}\\{#2}\end{pmatrix}}}
\newcommand{\dch}{\mathrm{dch}}

\newcommand{\zetaN}[3]{\zeta_{#1}{\begin{pmatrix}{#2}\\{#3}\end{pmatrix}}}

\newcommand{\zetam}{\zeta^{\mathfrak{m}}}

\newcommand{\zetamN}[3]{\zeta_{#1}^{\mathfrak{m}}{\begin{pmatrix}{#2}\\{#3}\end{pmatrix}}}
\newcommand{\zetaaN}[3]{\zeta_{#1}^{\mathfrak{a}}{\begin{pmatrix}{#2}\\{#3}\end{pmatrix}}}
\newcommand{\logm}{\log^{\mathfrak{m}}}
\renewcommand{\Im}{\mathrm{I}^{\mathfrak{m}}}
\newcommand{\Ia}{\mathrm{I}^{\mathfrak{a}}}

\newcommand{\myzetam}[2]{\zeta_{F_2}^{\mathfrak{m}}{\begin{pmatrix}{#1}\\{#2}\end{pmatrix}}}

\newcommand{\per}{\mathrm{per}}
\newcommand{\MZV}{\mathrm{MZV}}
\newcommand{\Gal}{\mathrm{Gal}}
\newcommand{\len}{\mathrm{len}}
\newcommand{\dep}{\mathrm{dep}}
\newcommand{\od}{\mathrm{odd}}
\newcommand{\ev}{\mathrm{even}}

\newcommand{\shu}{\hspace{-0.1mm}\text{\foreignlanguage{russian}{ш}}}

\DeclareFontFamily{U}{wncy}{}
\DeclareFontShape{U}{wncy}{m}{n}{<->wncyr10}{}
\DeclareSymbolFont{mcy}{U}{wncy}{m}{n}
\DeclareMathSymbol{\Shu}{\mathord}{mcy}{"58}

\newcommand{\bk}{\mathbf{k}}
\newcommand{\bp}{\boldsymbol{\varphi}}
\newcommand{\be}{\boldsymbol{\epsilon}}
\newcommand{\bla}{\boldsymbol{\lambda}}
\newcommand{\bl}{\mathbf{l}}
\newcommand{\br}{\mathbf{r}}

\newcommand{\id}{\mathrm{id}}
\newcommand{\Ker}{\mathrm{Ker}}

\address{Eisuke Otsuka\\
Mathematical inst. Tohoku Univ.\\
6-3, Aoba, Aramaki, Aoba-ku, Sendai, 980-8578, JAPAN}
\email{eisuke.otsuka.p3@dc.tohoku.ac.jp}

\begin{document}

\title{
On arithmetic properties of periods for some rational differential forms over $\bQ$ on the Fermat curve $F_2$ of degree 2
}
\author{
Eisuke Otsuka
}

\begin{abstract}
In this paper, we will define analogues of multiple zeta values by replacing the differential forms defining multiple zeta values with some $\bQ$-rational differential forms on the Fermat curve $F_2$ of degree 2 and discuss their arithmetic properties. We also investigate a motivic structure of the motivic periods corresponding to our periods. However, in order to study them, the current theory for motivic zeta elements is insufficient, and it leads us to study the base extension of the space of the motivic periods $\mathcal{H}_4$ of level 4 and its Galois invariant part.
\end{abstract}

\maketitle
\tableofcontents

\section{Introduction}
The multiple zeta values (MZVs for short, also called the Euler-Zagier type MZVs) are the real numbers defined by 
\begin{align*}
\zeta(k_1,\dots,k_d)=\sum_{0<n_1<\dots<n_d}\frac{1}{n_1^{k_1}\dots n_d^{k_d}}
\end{align*}
for $k_1,\dots,k_{d-1}\in\bZge{1}$, and $k_d\in\bZge{2}$. A pair $\bk=(k_1,\dots,k_d)\in(\bZge{1})^d$ with $k_d>1$ is called the (admissible) index of the weight $|\bk|:=k_1+\dots+k_d$. It is known that there are many $\bQ$-linear relations between MZVs (\hspace{-0.1mm}\cite{IKZ},\cite{LM},\cite{O} et al.), and in \cite{Z}, Zagier conjectured 
\begin{align} \label{eq: dim}
\sum_{k\ge0}d_kt^k=\frac{1}{1-t^2-t^3},
\end{align}
where $d_k$ is the dimension of the $\bQ$-linear space spanned by MZVs of weight $k$. In addition, MZVs appear in the study of quantum groups or knot invariants (\hspace{-0.1mm}\cite{Dr},\cite{LM}), and their arithmetic properties are now being actively studied (\hspace{-0.1mm}\cite{HS},\cite{Zh} et al.). One of the most important properties for MZVs is that they can be written as iterated integrals (\hspace{-0.1mm}\cite[p. 510]{Z}):
\begin{align}
\zeta(k_1,\dots,k_d)=\int_{0<x_1<\dots<x_k<1}\prod_{j=1}^k\phi_j(x_j), \label{eq: EZint}
\end{align}
where $$\phi_j(x)=\begin{cases}
\omega_1(x):=\frac{dx}{1-x}, & j=k_1+\dots+k_{s-1}+1 \text{ for some } s=1,\dots,d,\\
\omega_0(x):=\frac{dx}{x}, & \text{otherwise}.
\end{cases}$$ Put $\omega_1(x)=\omega_1$, $\omega_0(x)=\omega_0$ for simplicity.

In this paper, we replace $\omega_1$ in (\ref{eq: EZint}) with $\omega:=\frac{\diff x}{\sqrt{1-x^2}}$, and define the analogues of MZVs by 
\begin{align*}
\myzeta{k_1,\dots,k_d}{\varphi_1\dots,\varphi_d}:=\int_0^1\varphi_1\omega_0^{k_1-1}\dots\varphi_d\omega_0^{k_d-1}
\end{align*}
for $k_j\in\bZge{1}, \varphi\in\{\omega_1,\omega\}~(j=1,\dots,d)$. We call them multiple zeta values with respect to $F_2$ ($F_2$-MZVs for short) since $\omega$ is a $\bQ$-rational differential form on the Fermat curve $F_2$ of degree 2. These values have the iterated sum representation involving the central binomial coefficients (see Proposition \ref{prop: sum-rep}). We also give a motivic interpretation of $F_2$-MZVs to study their arithmetic algebraic properties, and analyze a motivic structure of $F_2$-MZVs. However, the current theory of motivic zeta elements with levels given by Deligne, Goncharov, Brown et al. is insufficient in our purpose, and therefore, we need to extend the base of the space of the motivic periods $\cH_4$ of level 4 and consider its Galois invariant part. We now explain our main results as below. We refer Section \ref{ss: DGB} for the notations.

\subsection{Main Results}
Let $\cH_4$ be the space of motivic periods of level 4, and $\widetilde{\cH}_4:=\cH_4\otimes_\bQ\bQ(\sqrt{-1})$ be its base extension. The Galois group $\Gal(\bQ(\sqrt{-1})/\bQ)=\{1,\sigma\}$ acts on $\widetilde\cH_4$ diagonally. Then the motivic interpretation of $F_2$-MZVs ($F_2$-MMZVs for short) is given by elements of $\left(\widetilde\cH_4\right)^\sigma$ which is $\sigma$-invariant part of $\widetilde\cH_4$. In addition, $\widetilde\cH_4$ is isomorphic non-canonically to the graded Hopf structure $\widetilde\cU_4$. We fix an isomorphism $\widetilde\Phi: \widetilde\cH_4\rightarrow \widetilde\cU_4$ defined in Proposition \ref{prop: regPhi}. The first result is about the $\sigma$-invariant part of $\widetilde\cU_4$.
\thm \label{theo: main1}
For $k\in\bZge{0}$, we have
\begin{align}
\left(\widetilde\cU^{(k)}\right)^\sigma=\bigoplus_{r\ge0}\left\langle f_{j_1}\dots f_{j_r}\{(2\pi\sqrt{-1})^{\mathfrak{m}}\}^l\otimes \sqrt{-1}^{k-r}~\middle|~\begin{array}{c} j_1,\dots,j_r\in\bZge{1}, l\in\bZge{0}\\ j_1+\dots+j_r+l=k\end{array}\right\rangle_\bQ,
\end{align}
where $\widetilde\cU^{(k)}$ is the homogeneous part of $\widetilde\cU_4$ in degree $k$ and $\left(\widetilde\cU^{(k)}\right)^\sigma$ is its $\sigma$-invariant part
\begin{align*}
\left(\widetilde\cU^{(k)}\right)^\sigma:=\left\{u\in\widetilde\cU^{(k)}~\middle|~\sigma u=u\right\}.
\end{align*}
\ethm
As an immediate application, we have the followings:
\cor[Corollary \ref{cor: dimension}]
Let $\MZV_{F_2}^{(k)}$ be the $\bQ$-linear space spanned by $F_2$-MZVs of weight $k$. Then it holds that
\begin{align*}
\dim_\bQ \MZV_{F_2}^{(k)}\le2^k.
\end{align*}
\ecor

We also have an application for direct sum decomposition of the $\bQ$-linear space $\MZV^{\mathfrak{m},(k)}_{F_2}$ spaned by $F_2$-MMZVs of weight $k$.
\cor[Corollary \ref{cor: 2}]
For each positive integer $k\ge1$, the following decomposition holds:
\begin{align*}
\MZV^{\mathfrak{m},(k)}_{F_2}&=\MZV^{\mathfrak{m},(k)}_{F_2,\mathrm{odd}}\oplus\MZV^{\mathfrak{m},(k)}_{F_2,\mathrm{even}},
\end{align*}
where
\begin{align*}
\MZV^{\mathfrak{m},(k)}_{F_2,\mathrm{odd}}&:=\left\langle\myzetam{k_1,\dots,k_d}{\varphi_1,\dots,\varphi_d}~\middle|~\substack{k_1+\dots+k_d=k,\\ \len_\omega(\bp)\mathrm{: odd}}\right\rangle_\bQ,\\
\MZV^{\mathfrak{m},(k)}_{F_2,\mathrm{even}}&:=\left\langle\myzetam{k_1,\dots,k_d}{\varphi_1,\dots,\varphi_d}~\middle|~\substack{k_1+\dots+k_d=k,\\ \len_\omega(\bp)\mathrm{: even}}\right\rangle_\bQ,
\end{align*}
and $\len_\omega(\bp):=\#\{j~|~\phi_j=\omega\}$. Also, we have
\begin{align*}
\dim_\bQ\MZV^{\mathfrak{m},(k)}_{F_2,\mathrm{odd}}&\le2^{k-1},\\
\dim_\bQ\MZV^{\mathfrak{m},(k)}_{F_2,\mathrm{even}}&\le2^{k-1}.
\end{align*}
\ecor

The second result is about an explicit formula for $F_2$-MMZVs of $d=1$.
\thm \label{theo: main2}
For each positive integer $k\ge1$, we have
\begin{align} \label{eq: main2}
\myzetam{k}{\omega}=\frac{(2\pi\sqrt{-1})^{\mathfrak{m}}}{4}\sum_{\mathbf{l}\in\mathbb{I}_{k-1}}\frac{(\logm2)^{l_1}}{l_1!}\prod_{j=2}^{|\mathbf{l}|}\frac{1}{l_j!}\left\{\frac{1-2^{1-j}}{j}\zetam(j)\right\}^{l_j}\otimes\frac{1}{\sqrt{-1}},
\end{align}
where $\pi^{\mathfrak{m}}=(2\pi\sqrt{-1})^{\mathfrak{m}}\otimes\frac{1}{2\sqrt{-1}}\in\widetilde{\cH}_4$, $\logm2$ is the motivic interpretation of $\log 2$ defined in Definition \ref{def: motivicLog}, and $\zetam(j)$ is a motivic multiple zeta value defined in Definition \ref{def: motivic-MZV}. Here we put
\begin{align*}
\begin{array}{ll}
\bI_0:=\{(0)\},&\\
\mathbb{I}_{k-1}:=\left\{(l_1,l_2,\dots,l_{s})\in\bigsqcup_{m>0}(\bZge{0})^m~\middle|~\sum_{j\ge1}jl_j=k-1,~~l_s\ne0\right\}, & k\ge2
\end{array},
\end{align*}
and $|\mathbf{l}|:=s$ for $\mathbf{l}=(l_1,l_2\dots, l_s)\in\mathbb{I}_{k-1}$. In addition, if $|\mathbf{l}|=1$, put $$\prod_{j=2}^{|\mathbf{l}|}\frac{1}{l_j!}\left\{\frac{1-2^{1-j}}{j}\zetam(j)\right\}^{l_j}=1.$$
\ethm

Note that the terms of
\begin{align*}
\sum_{\mathbf{l}\in\mathbb{I}_{k-1}}\frac{(\logm2)^{l_1}}{l_1!}\prod_{j=2}^{|\mathbf{l}|}\frac{1}{l_j!}\left\{\frac{1-2^{1-j}}{j}\zetam(j)\right\}^{l_j},
\end{align*}
in the right hand side of the equation (\ref{eq: main2}) are the weight $k-1$ elements of $\cH_4$ and the whole weight become $k$ together with the inisial $(2\pi\sqrt{-1})^\mathfrak{m}$.

In the course of the proof of Theorem \ref{theo: main2}, we first have an explicit formula of $F_2$-MZVs;
\begin{align} \label{eq: explicit_real}
\myzeta{k}{\omega}=\frac{\pi}{2}\sum_{\mathbf{l}\in\mathbb{I}_{k-1}}\frac{(\log2)^{l_1}}{l_1!}\prod_{j=2}^{|\mathbf{l}|}\frac{1}{l_j!}\left\{\frac{1-2^{1-j}}{j}\zeta(j)\right\}^{l_j}.
\end{align}
Then, we deduce Theorem \ref{theo: main2} by using (\ref{eq: explicit_real}) and our motivic structure of $F_2$-MZVs. As an immediate conseqence, since $\zeta(2n)\in\bQ\cdot\pi^{2n}$ for $n\in\bZge{1}$, the formula (\ref{eq: explicit_real}) yields
\begin{align*}
\myzeta{k}{\omega}\in\left\langle\pi^l(\log2)^{l_1}\prod_{j\ge3\text{: odd}}\zeta(j)^{l_j}~\middle|~\begin{array}{c}l,l_1,l_3,\dots\in\bZge{0}, l\text{: odd}\\ l+\sum_{j\ge1\text{: odd}}l_j=k\end{array}\right\rangle_\bQ.
\end{align*}

We organize this paper as follows. In Section \ref{ss: preparation}, we give some notations and definitions for $F_2$-MZVs. In Section \ref{ss: classical}, we discuss some basic properties of $F_2$-MZVs. Proposition \ref{prop: myzeta-level} among all is a key to study the motivic interpretaion of $F_2$-MZVs. Also, we give a proof of this formula (\ref{eq: explicit_real}) by using the iterated sum representation of $F_2$-MZVs in Proposition \ref{prop: explicit}. In Section \ref{ss: DGB}, we summarize the basic properties of motvic iterated integrals, and discuss the base extension and its Galois invariant part. We also discuss generalizations of path connection formula (Theorem \ref{theo: MII} (6)) and Goncharov's coaction formula (Theorem \ref{theo: coaction}) to prove the main theorems. After that, we begin to study $F_2$-MMZVs. In Section \ref{ss: proof}, we give proofs of the main results. Theorem \ref{theo: main1} is proved by calculating the Galois action on $\widetilde\cU_4$. We also calculate coaction of $\myzetam{k}{\omega}$ to prove Theorem \ref{theo: main2} by induction on $k$.

\section*{Acknowledgement}
In writing this paper, the author would like to express his deepest gratitude to his supervisor, Professor Takuya Yamauchi of Tohoku University. He gave the author a lot not only the basics of the geometric background of iterated integrals, but also many corrections to this paper. Though the author had no background in geometry, but Professor Yamauchi gave him very careful guidance and research problems fitting into his level. This enabled him to complete this paper. He would like to express his sincere gratitude to him.

He would like to thank Professor Yasuo Ohno for his guidance since last year, and for teaching him the fascination of multiple zeta values. He has participated in the seminar this year as well, and advised the author on various aspects of his research from various viewpoints.

He would also like to express his sincere gratitude to the members of Yamauchi Laboratory and Ohno Laboratory for their kindness both in his private and academic life. In particular, Dr. Yuya Murakami and Yuta Kadono gave him many advices on writing this paper. He is deeply grateful to them.

Finally, he would like to thank his family for many supports.

\section{Preparations} \label{ss: preparation}
In this section, we define some notations using in this paper, and give the definition of $F_2$-MZVs.
\nota
Through this paper, we use the following notation unless otherwise mentioned.
\begin{itemize}
  \item For a finite set $S$, let $\#S$ be the cardinality of $S$;
  \item Kronecker delta: for $r,s\in\bZ$, put $\delta_{r,s}:=\begin{cases}
  1, & r=s,\\
  0, & r\ne s;
  \end{cases}$
  \item binomial coefficients: for $k,r\in\bZge{0}$, put $\binom{k}{r}=\begin{cases}
  \frac{k!}{r!(k-r)!}, & k\ge r,\\
  0, & k<r;\end{cases}$
  \item Bernoulli numbers: for $n\ge0$, we define $B_n\in\bQ$ by the exponential generating function $$\sum_{n\ge0}\frac{B_n}{n!}x^n=\frac{x}{e^x-1};$$
  \item for $n,r,s\in\bZge{1}$, let $S_n$ be the $n$-th symmetry group, and define $(r,s)$-shuffle $S_{r,s}$ by $$S_{r,s}:=\left\{\delta\in S_{r+s}~\middle|~\substack{\delta^{-1}(1)<\dots<\delta^{-1}(r)\\ \delta^{-1}(r+1)<\dots<\delta^{-1}(r+s)}\right\};$$
  \item for $N\in\bZge{1}$, put $\xi_N:=\exp\left(\frac{2\pi\sqrt{-1}}{N}\right)$ and $i:=\xi_4=\sqrt{-1}$ for $N=4$;
  \item for $N\in\bZge{1}$, put $\mu_N:=\left\{\xi_N^r~\middle|~r=0,1,\dots,N-1\right\}$ and $\widetilde{\mu}_N:=\mu_N\cup\{0\}$;
  \item for $N\in\bZge{1}$, put $\bQ_N:=\bQ(\xi_N)$ and let $\mathcal{O}_{N}$ be its ring of integers;
  \item let $F_2:=\left\{[X:Y:Z]\in\bP^2(\bC)~\middle|~X^2+Y^2=Z^2\right\}$ be the Fermat curve of degree 2, and we regard it as an algebraic curve over $\bQ$ and a smooth scheme over $\bZ[1/2]$;
  \item let $x=X/Z, y=Y/Z$ be local variables of $F_2$;
  \item path reversal: for a smooth path $\gamma:[0,1]\rightarrow\bC$, let $\gamma^{-1}:[0,1]\rightarrow\bC$ be $$\gamma^{-1}(t):=\gamma(1-t);$$
  \item path connection: for smooth pathes $\gamma_1, \gamma_2:[0,1]\rightarrow\bC$ satisfying $\gamma_1(1)=\gamma_2(0)$, let $\gamma_1\gamma_2: [0,1]\rightarrow\bC$ be $$(\gamma_1\gamma_2)(t):=\begin{cases}
  \gamma_1(2t), & 0<t<1/2,\\
  \gamma_2(2t-1), & 1/2<t<1;
  \end{cases}$$
  \item we define the smooth path $\dch_{p,q}$ on $\bC$ by $\dch_{p,q}: [0,1]\rightarrow\bC, \dch_{p,q}(t)=p+t(q-p)$ for $p,q\in\bC$. In particular, for $p=0$ and $q=1$, we put $\dch:=\dch_{0,1}$.
\end{itemize}
We also use the empty sum in the summation symbol $\sum$ for $0$, and the empty product in the production symbol $\prod$ for $1$.
\enota

\defi[Iterated Integrals]
Let $M$ be a smooth manifold, $\gamma:[0,1]\rightarrow M$ be a smooth path on $M$, and $\phi_1, \dots, \phi_d$ be smooth $1$-forms on $M$. We define the iterated integral of $\phi_1, \dots, \phi_d$ along $\gamma$ by
\begin{align*}
\int_\gamma \phi_1\dots\phi_k:=\int_{0<t_1<\dots<t_k<1}\gamma^*\phi_1(t_1)\dots\gamma^*\phi_k(t_k).
\end{align*}
In particular, for $M=\bC$, $\gamma:[0,1]\rightarrow\bC$, and $a_1,\dots,a_k\in\bC\backslash\gamma((0,1))$ with $p:=\gamma(0)\ne a_1, q:=\gamma(1)\ne a_k$, we define
\begin{align} \label{def: II}
\mathrm{I}_\gamma(p;a_1,\dots,a_k;q):=\int_\gamma\frac{\diff x}{x-a_1}\dots\frac{\diff x}{x-a_k}.
\end{align}
In addition, for an element $w=e_{z_1}\dots e_{z_k}$ of the non-commutative polynomial ring $W_\bC:=\bQ\langle e_z~|~z\in\bC\rangle$ with $z_1\ne p, z_k\ne q$, we define
\begin{align*}
\mathrm{I}_\gamma(p;w;q):=\mathrm{I}_\gamma(p,z_1,\dots,z_k;q)
\end{align*}
and extend $\mathrm{I}_\gamma(p;*;q)$ to whole $W_\bC$ by the $\bQ$-linearity.
\edefi

\rem
By (\ref{eq: EZint}), the classical each Euler-Zagier type MZV has the following iterated integral representation
\begin{align*}
\zeta(k_1,\dots,k_d)=(-1)^d\mathrm{I}_\dch(0;1,\{0\}^{k_1-1},1,\{0\}^{k_2-1},\dots,1,\{0\}^{k_d-1};1),
\end{align*}
where $\{0\}^k=\overbrace{0,\cdots,0}^k$.
\erem

Now, we extend the definition of iterated integrals to the case when $p=a_1$ or $q=a_k$, and the integral in (\ref{def: II}) does not converge. 

\defi[Tangential Base Points]
A pair $(p,v)$ of a point $p\in\bC$ and a tangent vector $v\in T_p\bC\backslash\{0\}$ is called the tangential base point. A smooth path $\gamma$ from one tangential base point $(p,u)$ to another tangential base point $(q,v)$ is a smooth path $\gamma:[0,1]\rightarrow\bC$ satisfying
$$\gamma(0)=p,~~~ \gamma(1)=q,~~~ \gamma'(0)=u,~~~ \gamma'(1)=-v.$$
\edefi

\defi[Regularized Limit]
For a complex-valuesd function $f(\epsilon)$ defined for all sufficiently small real numbers $\epsilon>0$, we define
\begin{align*}
\underset{\epsilon\to0}{\mathrm{Reg}}~f(\epsilon):=c_{0,0}
\end{align*}
when $f(\epsilon)$ can be expressed as
\begin{align*}
f(\epsilon)=\sum_{s=0}^N\sum_{j=0}^\infty c_{s,j}(\log\epsilon)^s\epsilon^j
\end{align*}
for some $c_{s,j}\in\bC$. Note that this expression is unique exists.
\edefi

\defi[Generalization of Iterated Integral, {\cite[Definition 3.345]{GF}}] \label{def-regII}
For two tangential base points $\mathbf{p}=(p,u), \mathbf{q}=(q,v)$, a smooth path $\gamma:[0,1]\rightarrow\bC$ from $\mathbf{p}$ to $\mathbf{q}$, and $a_1,\dots,a_k\in\bC\backslash\gamma((0,1))$, we define the iterated integral by
\begin{align*}
\mathrm{I}_\gamma(\mathbf{p};a_1,\dots,a_k;\mathbf{q}):=\underset{\epsilon\to0}{\mathrm{Reg}}\int_{\gamma_\epsilon}\frac{\diff x}{x-a_1}\dots\frac{\diff x}{x-a_k},
\end{align*}
where $\gamma_\epsilon:=\gamma|_{[\delta(\epsilon),\delta'(\epsilon)]}$, $\delta(\epsilon):=\min\{t\in[0,1]~|~|p-\gamma(t)|\ge\epsilon\}$, and $\delta'(\epsilon):=\max\{t\in[0,1]~|~|q-\gamma(t)|\ge\epsilon\}$ for $0<\epsilon<1$. Clearly this is well-defined.
\edefi

Now, we define some rational 1-forms on $F_2$ to define $F_2$-MZVs.

\defi \label{def: phis}
Let $\omega_0=\frac{\diff x}{x},$ $\omega_1=\frac{\diff x}{1-x},$ and $\omega=\frac{\diff x}{y}$ be three rational $1$-forms on $F_2$. They are a part of the basis of the algebraic de Rham cohomology $H^1_{\mathrm{dR}}(Y_2)$ over $\bQ$. Here we put $Y_2:=\mathrm{Spec}\hspace{1mm} \bQ\left[x,y,\frac{1}{x(1-x)}\middle]\right/(x^2+y^2-1)$. A remaining basis is $\frac{x\diff x}{y}=x\omega$, but essentially there is no difference even if we use $x\omega$ below because $x\omega$ just shifts the weight by $-1$.
\edefi

We define $F_2$-MZVs as follows.

\defi[MZVs with respect to $F_2$ ($F_2$-MZVs for sohrt)] \label{def-myzeta}
Let $\gamma_1:[0,1]\rightarrow F_2(\bC)$ be a smooth path on $F_2(\bC)$ defined by $\gamma_1(t)=[t:\sqrt{1-t^2}:1]$. Then, for $\bk=(k_1,\dots,k_d)\in\bZge{1}^d$ and $\bp=(\varphi_1,\dots,\varphi_d)\in\{\omega_1,\omega\}^d$, we define $F_2$-MZVs by
\begin{align} \label{eq: myzeta}
\myzeta{k_1,\dots,k_d}{\varphi_1,\dots,\varphi_d}:=\int_{\gamma_1}\varphi_1\omega_0^{k_1-1}\dots\varphi_d\omega_0^{k_d-1}.
\end{align}
Here, $|\mathbf{k}|:=k_1+\dots+k_d$ is called the weight of $\bk$, $\dep(\bk):=d$ is called the depth of $\bk$, and $\len_\omega(\boldsymbol{\varphi}):=\#\{j~|~\varphi_j=\omega\}$ is called the length of $\bp$ respectively. Let $\MZV_{F_2}^{(k)}$ be the $\bQ$-linear space spanned by $F_2$-MZVs of weight $k$.
\edefi

\rem \label{rem: myzeta}
Regarding Definition \ref{def-myzeta}, we note that
\begin{enumerate}
\item the iterated integral (\ref{eq: myzeta}) converges if and only if $k_d>1$ or $\varphi_d=\omega$;
\item if $\varphi_j=\omega_1$ for all $j=1,\dots,d$, $\myzeta{k_1,\dots,k_d}{\omega_1,\dots,\omega_1}$ becomes the classical Euler-Zagier type MZVs $\zeta(k_1,\dots,k_d)$.
\end{enumerate}
\erem

By using Mathematica version 13.1, we can compute $F_2$-MZVs for $d=1$ and small $k$.
\exm
The following equations hold:
\begin{align*}
&\myzeta{1}{\omega}=\frac{\pi}{2},~~~
\myzeta{2}{\omega}=\frac{\pi}{2}\log2,~~~
\myzeta{3}{\omega}=\frac{\pi}{48}\left\{12(\log2)^2+\pi^2\right\},\\
&\myzeta{4}{\omega}=\frac{\pi}{48}\left\{4(\log2)^3+\pi^2\log2+6\zeta(3)\right\},\\
&\myzeta{5}{\omega}=\frac{\pi}{11520}\left\{240(\log2)^4+120\pi^2(\log2)^2+1440(\log2)\zeta(3)+19\pi^4\right\}.
\end{align*}
\eexm

One of our aims is to understand a reason why $\pi$, $\log2$, and odd zeta values show up in a complicated way.

\section{Classical Representations} \label{ss: classical}
In this section, we give some explicit representations of $F_2$-MZVs in terms of some analytic discussions.

\subsection{Another Expression of $\myzeta{k_1,\dots,k_d}{\varphi_1,\dots,\varphi_d}$}
\defi
Let $N\ge1$ be a positive integer, then for $k_1,\dots,k_d\in\bZge{1}$ and $\epsilon_1,\dots,\epsilon_d\in\mu_N$, the MZVs of level $N$ is defined by
\begin{align*}
\zetaN{N}{k_1,\dots,k_d}{\epsilon_1,\dots,\epsilon_d}&:=\sum_{0<n_1<\dots<n_d}\frac{\epsilon_1^{n_1}\dots\epsilon_d^{n_d}}{n_1^{k_1}\dots n_d^{k_d}}\\
&=(-1)^{d}\mathrm{I}_\dch(0;(\epsilon_1\dots\epsilon_d)^{-1},\{0\}^{k_1-1},\dots,\epsilon_d^{-1},\{0\}^{k_d-1};1).
\end{align*}
Let $\MZV_N^{(k)}$ be the $\bQ$-linear space spanned by MZVs of weight $k=k_1+\dots+k_d$.
\edefi

The following proposition shows that $F_2$-MZVs can be related to MZVs of level 4.
\propo \label{prop: myzeta-level}
Let $\bk\in(\bZge{1})^d$ and $\bp\in\{\omega_1,\omega\}^d$ with $(k_d,\varphi_d)\ne(1,\omega_1)$. Put $l=\len_\omega(\bp)$. Then, $F_2$-MZVs can be written as
\begin{align*}
\myzeta{k_1,\dots,k_d}{\varphi_1,\dots,\varphi_d}=(-\sqrt{-1})^l\mathrm{I}_\dch(0;\theta_1\eta_0^{k_1-1}\dots\theta_d\eta_0^{k_d-1};1),
\end{align*}
where
\begin{align*}
\theta_j=\begin{cases}
\eta:=e_i-e_{-i}, & \varphi_j=\omega,\\
\eta_1:=-2e_1+e_i+e_{-i}, & \varphi_j=\omega_1,
\end{cases}
\end{align*}
and $\eta_0=e_0-e_i-e_{-i}$. Furthermore, for weight $k$ index $(k_1,\dots,k_d)$,
\begin{align*}
\myzeta{k_1,\dots,k_d}{\varphi_1,\dots,\varphi_d}\in\sqrt{-1}^l\cdot\MZV_4^{(k)}.
\end{align*}
\epropo

\pf
Changing the variable as $t=\frac{2s}{1+s^2}$ in the definition of $F_2$-MZVs (the equation (\ref{eq: myzeta})), we have
\begin{align*}
\omega_0&=\frac{\diff t}{t}=\frac{(1-s^2)\diff s}{s(1+s^2)}=\frac{\diff s}{s}-\frac{\diff s}{s-\sqrt{-1}}-\frac{\diff s}{s+\sqrt{-1}},\\
\omega_1&=\frac{\diff t}{1-t}=\frac{2(1+s)\diff s}{(1-s)(1+s^2)}=-\frac{2\diff s}{s-1}+\frac{\diff s}{s-\sqrt{-1}}+\frac{\diff s}{s+\sqrt{-1}},\\
\omega&=\frac{\diff t}{\sqrt{1-t^2}}=\frac{2\diff s}{(1+s^2)}=-\sqrt{-1}\left(\frac{\diff s}{s-\sqrt{-1}}-\frac{\diff s}{s+\sqrt{-1}}\right),
\end{align*}
and the integral region is changed from $0<t_1<\dots<t_d<1$ to $0<s_1<\dots<s_d<1$. Thus, we have
\begin{align*}
\myzeta{k_1,\dots,k_d}{\varphi_1,\dots,\varphi_d}=(-\sqrt{-1})^l\mathrm{I}_\dch(0;\theta_1\eta_0^{k_1-1}\dots\theta_d\eta_0^{k_d-1};1).
\end{align*}
Furthermore, since $\mathrm{I}_\dch(0;\theta_1\eta_0^{k_1-1}\dots\theta_d\eta_0^{k_d-1};1)\in\MZV_4^{(k)}$, the above equation yields
\begin{align*}
\myzeta{k_1,\dots,k_d}{\varphi_1,\dots,\varphi_d}\in\sqrt{-1}^l\cdot\MZV_4^{(k)}.
\end{align*}
\epf

Now, we give the iterated sum representation of $F_2$-MZVs.

\propo \label{prop: sum-rep}
For $(k_1,\dots,k_d)\in\bZge{1}^d, (\varphi_1,\dots,\varphi_d)\in\{\omega_1,\omega\}^d$ with the convergence condition (Remark \ref{rem: myzeta} (1)), we have
\begin{align*}
\myzeta{k_1,\dots,k_d}{\varphi_1,\dots,\varphi_d}=\sum_{\substack{0<n_1<\dots<n_d\\n_1,\dots,n_d\in\bZge{1}}}\prod_{j=1}^d\frac{c(n_j-n_{j-1},\varphi_j)}{n_j^{k_j}},
\end{align*}
where $n_0:=0$ and the coefficients $c(n,\varphi)$ are defined by
\begin{align*} 
c(n,\varphi)=\begin{cases}
\frac{1}{2^{n-1}}\binom{n-1}{\frac{n-1}{2}}, & \varphi=\omega \text{ and } n\text{: odd},\\
0, & \varphi=\omega \text{ and } n\text{: even},\\
1, & \varphi=\omega_1
\end{cases}
\end{align*}
for $n\in\bZge{1}$ and $\varphi\in\{\omega_1,\omega\}$. In particular, in the case of $d=1$ and $\varphi_1=\omega$, we have
\begin{align} \label{eq: sumrep-d=1}
\myzeta{k}{\omega}=\sum_{n\ge 0}\frac{\binom{2n}{n}}{4^n(2n+1)^k}=1+\frac{1}{2}\cdot\frac{1}{3^k}+\frac{3}{8}\cdot\frac{1}{5^k}+\frac{5}{16}\cdot\frac{1}{7^k}+\cdots.
\end{align}
\epropo

\pf
Let $k=k_1+\dots+k_d$, $c(n,\omega_0):=\delta_{n,0}$ for a positive integer $n$, and $c(0,\omega_1)=c(0,\omega)=0$. For $\phi_j\in\{\omega_0,\omega_1,\omega\}$ ($j=1,\dots,k$) with $\phi_1\ne\omega_0$, $\phi_k\ne\omega_1$, and $\gamma_t=\gamma_1|_{[0,t]}$ ($0\le t\le1$), we will prove
\begin{align}
\int_{\gamma_t}\phi_1\dots\phi_k=\sum_{0< n_1\le\dots\le n_k}\prod_{j=1}^k\frac{c(n_j-n_{j-1},\phi_j)}{n_j}t^{n_k} \label{eq: sum-rep}
\end{align}
by induction on $k$.
Note that the Taylor series of $\omega_0=\frac{\diff x}{x}$, $\omega_1=\frac{\diff x}{1-x}$, and $\omega=\frac{\diff x}{y}=\frac{\diff x}{\sqrt{1-x^2}}$ at $x=0$ are given by
\begin{align*}
\frac{\diff x}{x}&=\sum_{n\ge0}\delta_{n,0}x^{n-1}\diff x=\sum_{n\ge0}c(n,\omega_0)x^{n-1}\diff x,\\
\frac{\diff x}{1-x}&=\sum_{n>0}x^{n-1}\diff x=\sum_{n\ge0}c(n,\omega_1)x^{n-1}\diff x,\\
\frac{\diff x}{\sqrt{1-x^2}}&=\sum_{m>0}\frac{1}{2^{2m}}\binom{2m}{m}x^{2m}\diff x=\sum_{n\ge0}c(n,\omega)x^{n-1}\diff x.
\end{align*}
When $k=1$, since $\phi_1\ne\omega_0, \omega_1$, and $\int_{\gamma_t}\phi_1=\int_{\gamma_t}\omega$ converges absolutely and uniformly, we can calculate the integration term by term, then
\begin{align*}
\int_{\gamma_t}\phi_1&=\int_{0<t_1<t}\sum_{n\ge0}c(n,\phi_1)t_1^{n-1}\diff t_1\\
&=\sum_{n\ge0}c(n,\phi_1)\int_{0<t_1<t}t_1^{n-1}\diff t_1\\
&=\sum_{n>0}\frac{c(n,\phi_1)}{n}t^n.
\end{align*}
When $k>1$, by the induction hypothesis,
\begin{align*}
\int_{\gamma_t}\phi_1\dots\phi_k
&=\int_{0<t_k<t}\left(\int_{0<t_1<\dots<t_{k-1}<t_k}\phi_1\dots\phi_{k-1}\right)\phi_k\\
&=\int_{0<t_k<t}\sum_{0< n_1\le\dots\le n_{k-1}}\prod_{j=1}^{k-1}\frac{c(n_j-n_{j-1},\phi_j)}{n_j}t_k^{n_{k-1}}\sum_{n\ge0}c(n,\phi_k)t_k^{n-1}\diff t_k\\
&=\sum_{0< n_1\le\dots\le n_{k-1}}\prod_{j=1}^{k-1}\frac{c(n_j-n_{j-1},\phi_j)}{n_j}t_k^{n_{k-1}}\sum_{n_k\ge n_{k-1}}c(n_k-n_{k-1},\phi_k)\int_{0<t_k<t}t_k^{n_k-1}\diff t_k\\
&=\sum_{0< n_1\le\dots\le n_k}\prod_{j=1}^k\frac{c(n_j-n_{j-1},\phi_j)}{n_j}t^{n_k},
\end{align*}
and we have the equation (\ref{eq: sum-rep}). Now, the clain follows by setting $t=1$ in the equation (\ref{eq: sum-rep}) because $c(n_j-n_{j-1},\omega_0)=\begin{cases}
1, & n_{j-1}=n_j,\\
0, & n_{j-1}<n_j,
\end{cases}$ and $c(n_j-n_{j-1},\phi)=0$ if $\phi\in\{\omega_1,\omega\}$ and $ n_{j-1}=n_j$.
\epf

\subsection{An Explicit Formula for the Depth One Case}
In this section, we give an explicit formula for $F_2$-MVZs of depth 1.

\propo \label{prop: explicit}
For each positive integer $k\ge1$, we have
\begin{align*}
\myzeta{k}{\omega}=\frac{\pi}{2}\sum_{\mathbf{l}\in\mathbb{I}_{k-1}}\frac{(\log2)^{l_1}}{l_1!}\prod_{j=2}^{|\mathbf{l}|}\frac{1}{l_j!}\left\{\frac{1-2^{1-j}}{j}\zeta(j)\right\}^{l_j}.
\end{align*}
where $\mathbb{I}_{k-1}$ is defined in Theorem \ref{theo: main2}.
\epropo

To prove this proposition, we prepare the following two lemmas.
\lem \label{lem: 3.5}
For each positive integer $k\ge1$, we have
\begin{align} \label{eq: zetabeta}
\myzeta{k}{\omega}=\frac{(-1)^{k-1}}{2^k(k-1)!}\frac{\diff^{k-1}}{\diff x^{k-1}}B\left(x,\frac{1}{2}\middle)\right|_{x=\frac{1}{2}},
\end{align}
where $B(x,y):=\int_0^1t^{x-1}(1-t)^{y-1}\diff t~~(x, y>0)$ is the beta function.
\elem

\pf
Calculating the integrals term by term, we have
\begin{align*}
\frac{(-1)^{k-1}}{2^k(k-1)!}\frac{\diff^{k-1}}{\diff x^{k-1}}B\left(x,\frac{1}{2}\right)
&=\frac{(-1)^{k-1}}{2^k(k-1)!}\frac{\diff^{k-1}}{\diff x^{k-1}}\int_0^1t^{x-1}(1-t)^{-\frac{1}{2}}\diff t\\
&=\frac{(-1)^{k-1}}{2^k(k-1)!}\frac{\diff^{k-1}}{\diff x^{k-1}}\int_0^1\sum_{n\ge0}\frac{\binom{2n}{n}}{4^n}t^{n+x-1}\diff t\\
&=\frac{(-1)^{k-1}}{2^k(k-1)!}\sum_{n\ge0}\frac{\binom{2n}{n}}{4^n}\frac{(-1)^{k-1}(k-1)!}{(n+x)^k}.
\end{align*}
Then the right hand side of the equation (\ref{eq: zetabeta}) is
\begin{align*}
\frac{(-1)^{k-1}}{2^k(k-1)!}\frac{\diff^{k-1}}{\diff x^{k-1}}B\left(x,\frac{1}{2}\middle)\right|_{x=\frac{1}{2}}
&=\frac{1}{2^k}\sum_{n\ge0}\frac{\binom{2n}{n}}{4^n}\frac{1}{(n+\frac{1}{2})^k}\\
&=\sum_{n\ge 0}\frac{\binom{2n}{n}}{4^n(2n+1)^k},
\end{align*}
which is as same as the iterated sum representation in Proposition \ref{prop: sum-rep} (the equation (\ref{eq: sumrep-d=1})).
\epf

\lem \label{lem: 3.6}
For each positive integer $k\ge1$, we have
\begin{align*}
\frac{\diff^{k-1}}{\diff x^{k-1}}B\left(x,\frac{1}{2}\right)=B\left(x,\frac{1}{2}\right)\sum_{\mathbf{r}\in\mathbb{I}_{k-1}}c_{\mathbf{r}}\prod_{j=1}^{|\mathbf{r}|}\left(\psi^{(j-1)}(x)-\psi^{(j-1)}\left(x+\frac{1}{2}\right)\right)^{r_j},
\end{align*}
where $\psi^{(j)}(x)$ is the polygamma function, which is defined by the logarithmic derivatives of the gamma function $\Gamma(x)=\int_0^\infty t^{x-1}e^{-t}\diff t$:
\begin{align*}
\psi(x)=\psi^{(0)}(x)=\frac{\Gamma'(x)}{\Gamma(x)},~~~\psi^{(j)}(x)=\frac{\diff^j}{\diff x^j}\psi^{(0)}(x).
\end{align*}
for $x>0$. Here, the coefficient $c_{\mathbf{r}}$ is defined by
\begin{align} \label{eq: coeff}
c_{\mathbf{r}}=\prod_{j=1}^{|\mathbf{r}|}\frac{(k-1)!}{(j!)^{r_j}r_j!}
\end{align}
for $\br\in\mathbb{I}_{k-1}$.
\elem

\pf
We give a proof by induction on $k$. When $k=1$, we have the trivial equation since $\bI_0=\{(0)\}$ and $c_0=1$. When $k>1$, put $F_j(x):=\psi^{(j)}(x)-\psi^{(j)}\left(x+\frac{1}{2}\right)$ for $j\in\bZge{1}$. Then, by the induction hypothesis, we have
\begin{align*}
\frac{\diff^k}{\diff x^k}B\left(x,\frac{1}{2}\right)&=\frac{\diff}{\diff x}B\left(x,\frac{1}{2}\right)\sum_{\mathbf{l}\in\mathbb{I}_{k-1}}c_{\mathbf{l}}\prod_{j=1}^{|\mathbf{l}|}F_{j-1}(x)^{l_j}\\
&=B\left(x,\frac{1}{2}\right)F_0(x)\sum_{\mathbf{l}\in\mathbb{I}_{k-1}}c_{\mathbf{l}}\prod_{j=1}^{|\mathbf{l}|}F_{j-1}(x)^{l_j}\\
&~~+B\left(x,\frac{1}{2}\right)\sum_{\mathbf{l}\in\mathbb{I}_{k-1}}c_{\mathbf{k}}\sum_{j=2}^{|\mathbf{l}|}F_0(x)^{l_1}\dots\frac{\diff}{\diff x}F_{j-2}(x)^{l_{j-1}}\dots F_{s-1}^{l_s}\\
&=B\left(x,\frac{1}{2}\right)F_0(x)\sum_{\mathbf{l}\in\mathbb{I}_{k-1}}c_{\mathbf{l}}\prod_{j=1}^{|\mathbf{l}|}F_{j-1}(x)^{l_j}\\
&~~+B\left(x,\frac{1}{2}\right)\sum_{\mathbf{l}\in\mathbb{I}_{k-1}}c_{\mathbf{l}}\sum_{j=2}^{|\mathbf{l}|}F_0(x)^{l_1}\dots l_{j-1}F_{j-2}(x)^{l_{j-1}-1}F_{j-1}(x)^{l_j+1}\dots F_s(x)^{l_s}.
\end{align*}
Now, we change the order of the sum: for each $\mathbf{r}\in\mathbb{I}_k$, we first calculate the sum of the terms $\bl\in\bI_{k-1}$ with $\br=(l_1+1,l_2,\dots,l_s)$ or $\br=(l_1,\dots,l_{j-1}-1,l_j+1,\dots,l_s)$ for some $j=2,\dots,|\bl|$. Then there exists $\tilde{c}_{\mathbf{r}}\in\bQ$ such that
\begin{align*}
\frac{\diff^k}{\diff x^k}B\left(x,\frac{1}{2}\right)=B\left(x,\frac{1}{2}\right)\sum_{\mathbf{r}\in\mathbb{I}_{k}}\tilde{c}_{\mathbf{r}}\prod_{j=1}^{|\mathbf{r}|}\left(F_{j-1}(x)\right)^{r_j},
\end{align*}
and the coefficient $\tilde{c}_{\mathbf{r}}$ satisfies
\begin{align*}
\tilde{c}_{r_1,\dots,r_s}=c_{r_1-1,\dots,r_s}+\sum_{m=2}^s(r_{m-1}+1)c_{r_1,\dots,r_{m-1}+1,r_m-1,\dots r_s}.
\end{align*}
Now, calculating the equation above and the equation (\ref{eq: coeff}), we have
\begin{align*}
\tilde{c}_{r_1,\dots,r_s}
&=c_{r_1-1,\dots,r_s}+\sum_{m=2}^s(r_{m-1}+1)c_{r_1,\dots,r_{m-1}+1,r_m-1,\dots r_s}\\
&=\prod_{j=1}^{|\mathbf{r}|}\frac{(k-1)!}{(j!)^{r_j}r_j!}\cdot r_1+\sum_{m=2}^s(r_{m-1}+1)\prod_{j=1}^{|\mathbf{r}|}\frac{(k-1)!}{(j!)^{r_j}r_j!}\cdot\frac{m!r_m}{(m-1)!(r_{m-1}+1)}\\
&=\prod_{j=1}^{|\mathbf{r}|}\frac{(k-1)!}{(j!)^{r_j}r_j!}\cdot(r_1+2r_2+\dots+ sr_s)\\
&=c_{r_1,\dots,r_s},
\end{align*}
and we have the claim.
\epf

\pf[of Proposition \ref{prop: explicit}]
By Lemma \ref{lem: 3.5}, \ref{lem: 3.6}, and $B\left(\frac{1}{2},\frac{1}{2}\right)=\pi$, the claim follows if we prove
\begin{align} \label{eq: pf321}
\begin{split}
\psi\left(\frac{1}{2}\right)-\psi(1)&=-2\log2\\
\psi^{(j-1)}\left(\frac{1}{2}\right)-\psi^{(j-1)}(1)&=(-1)^j(j-1)!(2^j-2)\zeta(j)
\end{split}
\end{align}
for $j>1$. First, by Legendre duplication formula for gamma function
\begin{align*}
\Gamma(2x)=\frac{2^{2x-1}}{\sqrt{\pi}}\Gamma(x)\Gamma\left(x+\frac{1}{2}\right),
\end{align*}
and its logarithmic derivative at $x=\frac{1}{2}$, we have
\begin{align*}
&2\psi(1)=2\log2+\psi\left(\frac{1}{2}\right)+\psi(1),\\
&\psi\left(\frac{1}{2}\right)-\psi(1)=-2\log2.
\end{align*}
Next, by series representation of gamma function
\begin{align*}
\psi^{(j-1)}(z)=(-1)^j(j-1)!\sum_{n=0}^\infty\frac{1}{(z+n)^j},
\end{align*}
for $j>1$, we have
\begin{align*}
\psi^{(j-1)}(1)&=(-1)^j(j-1)!\sum_{n=0}^\infty\frac{1}{(1+n)^j}=(-1)^j(j-1)!\zeta(j)\\
\psi^{(j-1)}\left(\frac{1}{2}\right)&=(-1)^j(j-1)!\sum_{n=0}^\infty\frac{1}{(\frac{1}{2}+n)^j}\\
&=(-1)^j(j-1)!2^j\sum_{n=0}^\infty\frac{1}{(2n+1)^j}\\
&=(-1)^j(j-1)!2^j\left(\sum_{n=1}^\infty\frac{1}{n^j}-\sum_{n=1}^\infty\frac{1}{(2n)^j}\right)\\
&=(-1)^j(j-1)!(2^j-1)\zeta(j),
\end{align*}
that yields the equation (\ref{eq: pf321}).
\epf

\section{Motivic Background} \label{ss: DGB}
In this section, we review some properties of motivic periods and motivic iterated integrals to study the geometric background of $F_2$-MZVs. Then, we consider base extension and Galois action to apply these theories to the motivic interpretation of $F_2$-MZVs.

\subsection{Motivic Iterated Integrals}
Here, we review some properties of motivic iterated integrals studied by Brown, Deligne, Glanois, Goncharov, Terasoma, and others. The next theorem is given by the theory of the mixed Tate motives over $\mathcal{O}_{N}\left[\frac{1}{N}\right]$.

\thm[\hspace{-0.1mm}{\cite[p. 211, 1.1]{Go}}]
For each positive integer $N\in\bZge{1}$, there exists a $\bQ$-algebra $\cH_N$ satisfying the following properties.

\begin{enumerate}
  \item The space $\cH_N=\bigoplus_{k\ge0}\cH_N^{(k)}$ is a graded $\bQ$-algebra, and $\cH_N^{(0)}=\bQ$.
  \item There exists a ring homomorphism $\per: \cH_N\rightarrow\bC$ called the period map.
  \item If $N>2$,
  \begin{enumerate}
    \item there exists an element $\tau=(2\pi\sqrt{-1})^\mathfrak{m}\in\cH_N^{(1)}$ such that $\per(\tau)=2\pi\sqrt{-1}$;
    \item the space $\cA_N:=\cH_N/(\tau\cdot\cH_N)$ has a structure of graded Hopf algebra. In particular, there exists a graded $\bQ$-linear homomorphism $$\Delta:\cA_N\rightarrow\cA_N\otimes_\bQ\cA_N$$ called the coproduct on $\cA_N$;
    \item the space $\cH_N$ has a structure of graded $\cA_N$-comodule. In particular, there exists a $\bQ$-linear map $$\Delta:\cH_N\rightarrow\cA_N\otimes_\bQ\cH_N$$ called the coaction on $\cH_N$.
  \end{enumerate}
  \item If $N=1, 2$,
\begin{enumerate}
    \item there exists an element $\tau^2=\{(2\pi\sqrt{-1})^\mathfrak{m}\}^2\in\cH_N^{(2)}$ such that $\per(\tau^2)=(2\pi\sqrt{-1})^2=-4\pi^2$;
    \item the space $\cA_N:=\cH_N/(\tau^2\cdot\cH_N)$ has a structure of graded Hopf algebra. In particular, there exists a graded $\bQ$-linear homomorphism $$\Delta:\cA_N\rightarrow\cA_N\otimes_\bQ\cA_N$$ called the coproduct on $\cA_N$;
    \item the space $\cH_N$ has a structure of graded $\cA_N$-comodule. In particular, there exists a $\bQ$-linear map $$\Delta:\cH_N\rightarrow\cA_N\otimes_\bQ\cH_N$$ called the coaction on $\cH_N$.
  \end{enumerate}
\end{enumerate}
Also, the coproduct $\Delta: \cA_N\rightarrow \cA_N\otimes_\bQ\cA_N$ is consistent with that induced by the coaction $\Delta: \cH_N\rightarrow\cA_N\otimes_\bQ\cH_N$.
Let $\rho: \cH_N\rightarrow\cA_N$ be the natural projection.
\ethm

\rem \label{rem: period_conjecture}
The period map $\per: \cH_N\rightarrow\bC$ is conjecturally injective (see \hspace{-0.1mm}{\cite[p.242, Conjecture 4.118]{GF}}), and the case $N=1$ includes Zagier's conjecture (the equation (\ref{eq: dim})).
\erem

The motivic iterated integrals are given as elements of $\cH_N$.
\thm[Motivic Iterated Integrals, {\cite[p.212, Theorem 1.1]{Go}}] \label{theo: MII}
For each positive integer $k\in\bZge{0}$, $k+2$ points $a_0,\dots,a_{k+1}\in\widetilde{\mu}_N$, two tangential base points $p=(a_0,u), q=(a_{k+1},v)$, and each path $\gamma$ from $p$ to $q$, there exists an element $$\Im_\gamma(p;a_1,\dots,a_k;q)\in\cH_N^{(k)}$$ called the motivic iterated integral, such that $$\per(\Im_\gamma(p;a_1,\dots,a_k;q))=\mathrm{I}_\gamma(p;a_1,\dots,a_k;q).$$ Furthermore, the motivic iterated integrals satisfy following conditions.
\begin{enumerate}
  \item The element $\rho(\Im_\gamma(p;a_1,\dots,a_k;q))$ is independent of the homotopy class of $\gamma$. We put $$\Ia(p;a_1,\dots,a_k;q):=\rho(\Im_\gamma(p;a_1,\dots,a_k;q))\in\cA_N^{(k)};$$
  \item if $k=0$, $\Im_\gamma(p;q)=1$;
  \item if $k>0$ and $p=q$, $\Im_\gamma(p;a_1,\dots,a_k;p)=0$;
  \item the equations $\Im_\dch(0;0;1)=\Im_\dch(0;1;1)=0$ hold;
  \item (path reversal): for a smooth path $\gamma$ on $\bC$, $\Im_{\gamma^{-1}}(p;a_1,\dots,a_k;q)=(-1)^k\Im_\gamma(q;a_k,\dots,a_1;p)$;
  \item (path connection): for smooth paths $\gamma_1, \gamma_2$ on $\bC$, \\
  $\Im_{\gamma_1\gamma_2}(p;a_1,\dots,a_k;q)=\sum_{s=1}^k\Im_{\gamma_1}(p;a_1,\dots,a_s;\gamma_1(1))\Im_{\gamma_2}(\gamma_2(0);a_{s+1},\dots,a_k;q)$;
  \item (substitution): for $c\in\mu_N$, $\Im_\gamma(p;a_1,\dots,a_k;q)=\Im_{c\gamma}(cp;ca_1,\dots,ca_k;cq)$, where $c\gamma: [0,1]\rightarrow \bC$ is defined by $(c\gamma)(t):=c\gamma(t)$.
\end{enumerate}
\ethm

The path connection formula (Theorem \ref{theo: MII} (6)) can be generalized as follows.
\propo \label{prop: pathconnect}
For $r,k\in\bZge{1}$ and smooth paths $\gamma_1,\dots,\gamma_r$ satisfying $\gamma_j(1)=\gamma_{j+1}(0)~(j=1,\dots,r-1)$, 
\begin{align*}
\Im_\gamma(\gamma(0);a_1,\dots,a_k;\gamma(1))=\sum_{0\le j_1\le \dots\le j_{r-1}\le k}\prod_{s=1}^r\Im_{\gamma_s}(\gamma_s(0);a_{j_{s-1}+1},\dots,a_{j_s};\gamma_s(1)),
\end{align*}
where $j_0=0$, $j_r=k$, and $\gamma=\gamma_1\dots\gamma_r$.
\epropo
\pf
We prove the claim by induction on $r$. For $r=1$ we have the trivial equation. For $r>1$, by path connection formula and the induction hypothesis, we have
\begin{align*}
&\Im_\gamma(\gamma(0);a_1,\dots,a_k;\gamma(1))\\
&=\sum_{j_{r-1}=0}^k\Im_{\gamma_1\dots\gamma_{r-1}}(\gamma(0);a_1,\dots a_{j_{r-1}};\gamma_{r-1}(1))\cdot\Im_{\gamma_r}(\gamma_r(0);a_{j_{r-1}+1},\dots a_{k};\gamma(1))\\
&=\sum_{j_{r-1}=0}^k\sum_{0\le j_1\le \dots\le j_{r-2}\le j_{r-1}}\left\{\prod_{s=1}^{r-1}\Im_{\gamma_s}(\gamma_s(0);a_{j_{s-1}+1},\dots,a_{j_s};\gamma_s(1))\right\}\cdot\Im_{\gamma_r}(\gamma_r(0);a_{j_{r-1}+1},\dots a_{k};\gamma(1))\\
&=\sum_{0\le j_1\le \dots\le j_{r-1}\le k}\prod_{s=1}^r\Im_{\gamma_s}(\gamma_s(0);a_{j_{s-1}+1},\dots,a_{j_s};\gamma_s(1)).
\end{align*}
then the claim follows.
\epf

The coaction $\Delta$ on $\cH_N$ for the motivic iterated integrals has an explicit formula given by Goncharov as follows.
\thm[Coaction Formula, {\cite[p. 213, Theorem 1.2]{Go}}] \label{theo: coaction}
The following formula holds:
\begin{align*}
\Delta\Im_\gamma(a_0;a_1,\dots,a_k;a_{k+1})&=\sum_{r=0}^k\sum_{0=j_0<j_1<\dots<j_r<j_{r+1}=k+1}\\
&\left(\prod_{s=0}^r\Ia(a_{j_s};a_{j_s+1},\dots,a_{j_{s+1}-1};a_{j_{s+1}})\right)\otimes\Im_\gamma(a_0;a_{j_1},\dots,a_{j_r};a_{k+1}).
\end{align*}
\ethm

\defi \label{def Deltas}
For the coaction $\Delta: \cH_N\rightarrow\cA_N\otimes_\bQ\cH_N$ on $\cH_N$, let $\Delta', \widetilde\Delta: \cH_N\rightarrow\cA_N\otimes_\bQ\cH_N$ be two $\bQ$-linear maps defined by
\begin{align*}
\Delta'&:=\Delta-1\otimes\id_{\cH_N},\\
\widetilde\Delta&:=\Delta-1\otimes\id_{\cH_N}-\rho\otimes1.
\end{align*}
We also define two $\bQ$-linear maps $\Delta', \widetilde\Delta: \cA_N\rightarrow\cA_N\otimes_\bQ\cA_N$ by
\begin{align*}
\Delta'&:=\Delta-1\otimes\id_{\cA_N},\\
\widetilde\Delta&:=\Delta-1\otimes\id_{\cA_N}-\id_{\cA_N}\otimes1.
\end{align*}
\edefi
Sometimes we don't calculate $\Delta$ itself, but $\Delta'$ or $\widetilde\Delta$.

\rem \label{rem: Deltas}
By Theorem \ref{theo: coaction}, $\Delta'\Im_\gamma(a_0;a_1,\dots,a_k;a_{k+1})$ and $\widetilde\Delta\Im_\gamma(a_0;a_1,\dots,a_k;a_{k+1})$ can be written as follows.
\begin{align*}
\Delta'\Im_\gamma(a_0;a_1,\dots,a_k;a_{k+1})&=\sum_{r=0}^{k-1}\sum_{0=j_0<j_1<\dots<j_r<j_{r+1}=k+1}\\
&\left(\prod_{s=0}^r\Ia(a_{j_s};a_{j_s+1},\dots,a_{j_{s+1}-1};a_{j_{s+1}})\right)\otimes\Im_\gamma(a_0;a_{j_1},\dots,a_{j_r};a_{k+1}),\\
\widetilde\Delta\Im_\gamma(a_0;a_1,\dots,a_k;a_{k+1})&=\sum_{r=1}^{k-1}\sum_{0=j_0<j_1<\dots<j_r<j_{r+1}=k+1}\\
&\left(\prod_{s=0}^r\Ia(a_{j_s};a_{j_s+1},\dots,a_{j_{s+1}-1};a_{j_{s+1}})\right)\otimes\Im_\gamma(a_0;a_{j_1},\dots,a_{j_r};a_{k+1}).
\end{align*}
\erem

The motivic multiple zeta values are defined by using motivic iterated integrals.
\defi[Motivic Multiple Zeta Values \cite{DG}, cf. {\cite[p. 344, Definition 2.2]{Gl}}] \label{def: motivic-MZV}
Let $N\ge1$ be a positive integer. Then, for $k_1,\dots,k_d\in\bZge{1}$ and $\epsilon_1,\dots,\epsilon_d\in\mu_N$, the motivic multiple zeta values (MMZVs for short) of level $N$ is defined by
\begin{align*}
\zetamN{N}{k_1,\dots,k_d}{\epsilon_1,\dots,\epsilon_d}:=(-1)^{d}\Im_\dch(0;(\epsilon_1\dots\epsilon_d)^{-1},\{0\}^{k_1-1},\dots,\epsilon_d^{-1},\{0\}^{k_d-1};1).
\end{align*}
In particular, when $N=1$, put
\begin{align*}
\zetam(k_1,\dots,k_d)&:=\zetamN{1}{k_1,\dots,k_d}{1,\dots,1}\\
&=(-1)^{d}\Im_\dch(0;(1,\{0\}^{k_1-1},\dots,1,\{0\}^{k_d-1};1)
\end{align*}
and
\begin{align*}
\mathcal{Z}_N^{(k)}:=\left\langle\zetamN{N}{k_1,\dots,k_d}{\epsilon_1,\dots,\epsilon_d}~\middle|~k_1+\dots+k_d=k\right\rangle_{\bQ}.
\end{align*}
By definition, we have $\per\left(\zetamN{N}{k_1,\dots,k_d}{\epsilon_1,\dots,\epsilon_d}\right)=\zetaN{N}{k_1,\dots,k_d}{\epsilon_1,\dots,\epsilon_d}$.
\edefi

MMZVs have some $\bQ$-linear relations as follows.
\propo[\hspace{-0.1mm}{\cite[p. 959, Lemma 3.2]{B}}] \label{prop: zetamev}
For each positive even number $k\ge1$, we have
\begin{align*}
\zetam(k)=-\frac{B_k}{2k!}\tau^k.
\end{align*}
\epropo

\propo[\hspace{-0.1mm}{\cite[Lemma 3.1]{Gl}}] \label{prop: conjugate}
For each positive integer $k\in\bZge{1}$ and $\epsilon\in\mu_N$, we have
\begin{align*}
\zetaaN{N}{k}{\epsilon^{-1}}=(-1)^{k-1}\zetaaN{N}{k}{\epsilon}.
\end{align*}
\epropo

\propo[\hspace{-0.1mm}{\cite[Lemma 3.1]{Gl}}] \label{prop: distribution}
For positive integers $n,N\ge1$ with $n|N$, put $M=\frac{N}{n}$. Then, for $\bk\in(\bZge{1})^d$ and $\be\in(\mu_{M})^d$, we have
\begin{align*}
\zetaaN{M}{\bk}{\be}=n^{|\bk|-d}\sum_{\bla}\zetaaN{N}{\bk}{\bla},
\end{align*}
where the index $\bla=(\lambda_1,\dots,\lambda_d)\in(\mu_N)^d$ runs with $\lambda_j^n=\epsilon_j~(j=1,\dots,d)$.
\epropo

The following properties hold for MMZVs of level $1$, $2$, and $4$. They are mainly used for the proofs of the main results.

\thm[\hspace{-0.1mm}{\cite[Theorem 8.1]{B}} for $N=1$, \cite{De} for $N=2,4$] \label{theo: BD}
If $N=1,2,4$, the following equation holds.
\begin{align*}
\mathcal{Z}_N=\cH_N.
\end{align*}
\ethm

\defi
For $N=1,2,4$ and $r\ge0$, put
\begin{align*}
\cD_r\cH_N^{(k)}:=\left\langle\zetamN{N}{k_1,\dots,k_d}{\epsilon_1,\dots,\epsilon_d}~\middle|~k_1+\dots+k_d=k, d\le r\right\rangle_{\bQ}.
\end{align*}
Then, $\cD=\{\cD_r\}_{r\ge0}$ becomes a filtration on $\cH_N$ by Theorem \ref{theo: BD}.
\edefi

\propo[\hspace{-0.1mm}{\cite[Corollary 3.7]{Gl}}] \label{prop: kernelM}
If $N=2,4$, for the coaction $\Delta$ on $\cH_N$, we have
\begin{align*}
\Ker\widetilde\Delta\cap\cH_2^{(k)}&=\bQ\zetamN{2}{k}{-1},\\
\Ker\widetilde\Delta\cap\cH_4^{(k)}&=\bQ\tau^k\oplus\bQ\zetamN{4}{k}{\sqrt{-1}}.
\end{align*}
\epropo

Now, we define the motivic logarithms, which are motivic analogues of logarithms.
\defi \label{def: motivicLog}
For each positive integer $N\in\bZge{1}$ and $\epsilon\in\mu_N$, we define $\logm(1-\epsilon)\in\cH_N^{(1)}$ by
\begin{align*}
\logm(1-\epsilon):=\Im_{\dch_{0,\epsilon}}(0;1;\epsilon).
\end{align*}
By definition, we have $\per(\logm(1-\epsilon))=\log(1-\epsilon)$.
\edefi

\exm
If $\epsilon=-1$ in Definition \ref{def: motivicLog}, we have
\begin{align*}
\logm2=\Im_{\dch_{0,-1}}(0;1;-1)=\Im_{\dch}(0;-1;1)=-\zetamN{2}{1}{-1}
\end{align*}
as an element of $\cH_2^{(1)}$.
\eexm

\subsection{Structure Theorem}
The Hopf strucute of $\cH_N$ has been determined. The next theorem is one example for the case of $N=1,2,4$.

\thm[Structure Theorem, {\cite[Lemma 2.1]{Gl}}] \label{theo: str}
If $N=1,2,4$, there exists a non-canonical isomorphism as non-commutative graded $\bQ$-algebras\begin{align*}
\begin{array}{rccc}
\Phi: & \cH_N&\xrightarrow{\sim}&\cU_N\\
&\rotatebox{-90}{\hspace{-1mm}$\xrightarrow[\rotatebox{90}{$\rho$}]{}$}& & \rotatebox{-90}{\hspace{-1mm}$\xrightarrow{\rotatebox{90}{$\rho$}}$}\\
\Phi: &\cA_N&\xrightarrow[\sim]{}&\cV_N,
\end{array}
\end{align*}
which are commutative with the projection $\rho: \cH_N\rightarrow\cA_N$. Here, $\cU_N$ and $\cV_N$ are the non-commutative polynomial rings
\begin{align*}
\cU_1&=\bQ\langle f_j~|~j\in\bZge{3}\text{: odd}\rangle\otimes_\bQ\bQ[\tau^2],\\
\cV_1&=\bQ\langle f_j~|~j\in\bZge{3}\text{: odd}\rangle,\\
\cU_2&=\bQ\langle f_j~|~j\in\bZge{1}\text{: odd}\rangle\otimes_\bQ\bQ[\tau^2],\\
\cV_2&=\bQ\langle f_j~|~j\in\bZge{1}\text{: odd}\rangle,\\
\cU_4&=\bQ\langle f_j~|~j\in\bZge{1}\rangle\otimes_\bQ\bQ[\tau],\\
\cV_4&=\bQ\langle f_j~|~j\in\bZge{1}\rangle,
\end{align*}
spanned by degree $j$ elements $f_j$, and $\rho: \cU_N\rightarrow\cV_N$ are defined by $$\rho(f_{j_1}\cdots f_{j_r}\tau^l)=\begin{cases}
f_{j_1}\cdots f_{j_r}, & l=0,\\
0, & l>0.
\end{cases}$$ The product $\shu: \cU_N\otimes_\bQ\cU_N\rightarrow\cU_N$ and the coaction $\Delta: \cU_N\rightarrow\cV_N\otimes_\bQ\cU_N$ on $\cU_N$ are given by
\begin{align}
\label{eq: shuU} &(f_{j_1}\dots f_{j_r}\tau^l)\shu(f_{j_{r+1}}\dots f_{j_{r+s}}\tau^m)=\sum_{\delta\in S_{r,s}}f_{j_{\delta(1)}}\dots f_{j_{\delta(r+s)}}\tau^{l+m};\\
\label{eq: coaU} &\Delta(f_{j_1}\dots f_{j_r}\tau^l)=\sum_{s=0}^rf_{j_1}\dots f_{j_s}\otimes f_{j_{s+1}}\dots f_{j_r}\tau^l.
\end{align}
Also, let $\Delta', \widetilde\Delta: \cU_N\rightarrow\cV_N\otimes_\bQ\cU_N$ be two $\bQ$-linear maps defined by
\begin{align*}
\Delta'&:=\Delta-1\otimes\id_{\cU_N},\\
\widetilde\Delta&:=\Delta-1\otimes\id_{\cU_N}-\rho\otimes1.
\end{align*}
\ethm

\lem \label{lemm: kernelH}
For the coaction $\Delta$ on $\cU_4$, we have
\begin{align*}
\Ker\Delta'\cap\cU_4^{(k)}&=\bQ \tau^k,\\
\Ker\widetilde\Delta\cap\cU_4^{(k)}&=\bQ \tau^k\oplus\bQ f_k.
\end{align*}
\elem
\pf
Since $\Delta$ is injective, it is sufficient if we investigate the case of $u=f_{j_1}\dots f_{j_r}\tau^l\in\cU_4^{(k)}$. Now, $\cU_N$ and $\cV_N$ are free $\bQ$-algebras, $\cV_N\otimes_\bQ\cU_N$ is also free. Then, each terms $f_{j_1}\dots f_{j_s}\otimes f_{j_{s+1}}\dots f_{j_r}\tau^l~~(s=1,\dots r)$, which appear in the definition of $\Delta'(u)$, are $\bQ$-linear independent on $\bigoplus_{m=0}^k(\cV_4^{(m)}\otimes_\bQ\cU_4^{(k-m)})$, then we have
\begin{align*}
\Delta'(u)=0 &\Longrightarrow \sum_{s=1}^rf_{j_1}\dots f_{j_s}\otimes f_{j_{s+1}}\dots f_{j_r}\tau^l=0\\
&\Longrightarrow r=0\\
&\Longrightarrow u=\tau^l=\tau^k.
\end{align*}
On the other hands, if $\widetilde\Delta(u)=0$, we have
\begin{align*}
\widetilde\Delta(u)=0 &\Longrightarrow \sum_{s=1}^rf_{j_1}\dots f_{j_s}\otimes f_{j_{s+1}}\dots f_{j_r}\tau^l=\begin{cases}
0, & l>0,\\
f_{j_1}\dots f_{j_r}\otimes1, & l=0,
\end{cases}\\
&\Longrightarrow r=0 \text{ or } l=0\\
&\Longrightarrow u=\tau^l=\tau^k \text{ or } u=f_{j_1}\dots f_{j_r}.
\end{align*}
\epf

\subsection{Extension of the Notations for Motivic Iterated Integrals}
For a subset $X\ne\emptyset$ of $\bC$, let $W_X:=\bQ\langle e_x~|~x\in X\rangle$ be the non-commutative polynomial ring spanned by the degree $1$ elements $e_x~(x\in X)$. Let $W_X^{(k)}$ be the degree $k$ homogeneous part of $W_X$.
\defi \label{def: shuffle}
For an element $w=e_{z_1}\dots e_{z_k}\in W_{\widetilde\mu_N}$, the motivic iterated integrals $\Im_\gamma(p;w;q)$ are defined by
\begin{align*}
\Im_\gamma(p;w;q):=\Im_\gamma(p;z_1,\dots,z_k;q),
\end{align*}
and extend to whole $W_{\widetilde\mu_N}$ by $\bQ$-linearlity.
\edefi

\propo[Shuffle Relation, \hspace{-0.1mm}{\cite{Go}}]
For two elements $w_1,w_2\in W_X$, the following equation holds:
\begin{align*}
\Im_\gamma(p;w_1;q)\Im_\gamma(p;w_2;q)=\Im_\gamma(p;w_1\shu w_2;q),
\end{align*}
where $\shu: W_X\otimes_\bQ W_X\rightarrow W_X$ is defined by
\begin{align} \label{def-shufflewords}
(e_{x_1}\dots e_{x_r})\shu(e_{x_{r+1}}\dots e_{x_{r+s}})=\sum_{\delta\in S_{r,s}}e_{x_{\delta(1)}}\dots e_{x_{\delta(r+s)}}
\end{align}
and the $\bQ$-linearlity of $W_X$.
\epropo

\lem \label{lemm: exp}
For each positive integer $k\in\bZge1$ and $w\in W_{\widetilde\mu_N}^{(1)}$, 
\begin{align*}
w^k=\frac{1}{k!}w^{\shu\hspace{0mm}k},
\end{align*}
where $w^k=\underbrace{w\cdots w}_{k}$ is the $k$ times non-commutative production of $w$, and $w^{\shu \hspace{0mm} k}=\underbrace{w\shu\cdots\shu w}_{k}$ is the $k$ times shuffle production of $w$.
\elem
\pf
We prove the claim by induction on $k$. When $k=1$, we have the trivial equation. When $k>1$, put $w=\sum_{z\in\widetilde\mu_N}c_ze_z$, then we have
\begin{align*}
w^{\shu\hspace{0mm}k}
=w\shu w^{\shu\hspace{0mm}k-1}
&=(k-1)!\cdot w\shu w^{k-1}\\
&=(k-1)! \sum_{z_1,\dots,z_k\in\widetilde\mu_N}\left(\prod_{s=1}^kc_{z_s}\right)(e_{z_1}\shu e_{z_2}\cdots e_{z_k})
\end{align*}
by the induction hypothesis. Now, by the definition of $\shu$ (the equation (\ref{def-shufflewords})),
$$e_{z_1}\shu e_{z_2}\cdots e_{z_k}=\sum_{s=1}^k e_{z_2}\cdots e_{z_s}e_{z_1}e_{z_{s+1}}\cdots e_{z_k}.$$ Therefore, we have
\begin{align*}
w^{\shu\hspace{0mm}k}
&=(k-1)! \sum_{z_1,\dots,z_k\in\widetilde\mu_N}\left(\prod_{s=1}^kc_{z_s}\right)(e_{z_1}\shu e_{z_2}\cdots e_{z_k})\\
&=(k-1)! \sum_{z_1,\dots,z_k\in\widetilde\mu_N}\left(\prod_{s=1}^kc_{z_s}\right)\left(\sum_{s=1}^k e_{z_2}\cdots e_{z_s}e_{z_1}e_{z_{s+1}}\cdots e_{z_k}\right)\\
&=k! \sum_{z_1,\dots,z_k\in\widetilde\mu_N}\left(\prod_{s=1}^kc_{z_s}\right)\left(e_{z_1}\cdots e_{z_k}\right)\\
&=k!\cdot w^k.
\end{align*}
\epf

The following proposition is the extension of Theorem \ref{theo: coaction} for the motivic iterated integrals on $W_{\widetilde{\mu}_N}$. We use this explicit formula in the proof of Theorem \ref{theo: main2}.
\propo \label{prop: coaction2}
Let $N\ge1$ be a positive integer, $c_z^{(j)}\in\bQ$, and $w_j=\sum_{z\in\widetilde\mu_N}c_z^{(j)}e_z\in W_{\widetilde\mu_N}^{(1)}$. Then for $w=w_1\cdots w_k\in W_{\widetilde\mu_N}$, we have
\begin{align*}
\Delta\Im_\gamma(p;w;q)&=\sum_{r=0}^k\sum_{(a_1,\dots,a_r)\in\widetilde\mu_N^r}\sum_{0=j_0<j_1<\dots<j_r<j_{r+1}=k+1}\\
&\left(\prod_{s=1}^rc_{a_s}^{(j_s)}\right)\cdot\left(\prod_{s=0}^r\Ia(a_s;w_{j_s+1}\cdots w_{j_{s+1}-1};a_{s+1})\right)\otimes\Im_\gamma(p;a_1,\dots,a_r;q),\\
\Delta\Ia(p;w;q)&=\sum_{r=0}^k\sum_{(a_1,\dots,a_r)\in\widetilde\mu_N^r}\sum_{0=j_0<j_1<\dots<j_r<j_{r+1}=k+1}\\
&\left(\prod_{s=1}^rc_{a_s}^{(j_s)}\right)\cdot\left(\prod_{s=0}^r\Ia(a_s;w_{j_s+1}\cdots w_{j_{s+1}-1};a_{s+1})\right)\otimes\Ia(p;a_1,\dots,a_r;q),
\end{align*}
where $a_0=p, a_{k+1}=q$. Moreover, we have
\begin{align*}
\widetilde\Delta\Im_\gamma(p;w;q)&=\sum_{r=1}^{k-1}\sum_{(a_1,\dots,a_r)\in\widetilde\mu_N^r}\sum_{0=j_0<j_1<\dots<j_r<j_{r+1}=k+1}\\
&\left(\prod_{s=1}^rc_{a_s}^{(j_s)}\right)\cdot\left(\prod_{s=0}^r\Ia(a_s;w_{j_s+1}\cdots w_{j_{s+1}-1};a_{s+1})\right)\otimes\Im_\gamma(p;a_1,\dots,a_r;q),\\\widetilde\Delta\Ia(p;w;q)&=\sum_{r=1}^{k-1}\sum_{(a_1,\dots,a_r)\in\widetilde\mu_N^r}\sum_{0=j_0<j_1<\dots<j_r<j_{r+1}=k+1}\\
&\left(\prod_{s=1}^rc_{a_s}^{(j_s)}\right)\cdot\left(\prod_{s=0}^r\Ia(a_s;w_{j_s+1}\cdots w_{j_{s+1}-1};a_{s+1})\right)\otimes\Ia(p;a_1,\dots,a_r;q).
\end{align*}
\epropo
\pf
By Theorem \ref{theo: coaction}, we have
\begin{align*}
\Delta\Im_\gamma(p;w;q)&=\Delta\Im_\gamma(p;w_1\dots w_k;q)\\
&=\Delta\sum_{\mathbf{z}\in\widetilde\mu_N^k}\left(\prod_{j=1}^kc_{z_j}^{(j)}\right)\Im_\gamma(p;\mathbf{z};q)\\
&=\sum_{\mathbf{z}\in\widetilde\mu_N^k}\left(\prod_{j=1}^kc_{z_j}^{(j)}\right)\sum_{r=0}^k\sum_{0=j_0<j_1<\dots<j_r<j_{r+1}=k+1}\\
&~~~\prod_{s=0}^r\Ia(z_{j_s};z_{j_s+1},\dots,z_{j_{s+1}-1};z_{j_{s+1}})\otimes\Im_\gamma(p;z_{j_1},\dots, z_{j_r};q)\\
&=\sum_{r=0}^k\sum_{\mathbf{z}\in\widetilde\mu_N^k}\sum_{0=j_0<j_1<\dots<j_r<j_{r+1}=k+1}\\
&~~~\left(\prod_{j=1}^kc_{z_j}^{(j)}\right)\cdot\prod_{s=0}^r\Ia(z_{j_s};z_{j_s+1},\dots,z_{j_{s+1}-1};z_{j_{s+1}})\otimes\Im_\gamma(p;z_{j_1},\dots, z_{j_r};q).
\end{align*}
Now, we change the order of the sum: for each $(a_1,\dots,a_r)\in\widetilde\mu_N^r$, we calculate the sum of the terms with $(z_{j_1},\dots,z_{j_r})=(a_1,\dots,a_r)$, then
\begin{align*}
\Delta\sum_{\mathbf{z}\in\widetilde\mu_N^k}\prod_{j=1}^kc_{z_j}^{(j)}\Im_\gamma(p;\mathbf{z};q)&=\sum_{r=0}^k\sum_{(a_1,\dots,a_r)\in\widetilde\mu_N^r}\sum_{0=j_0<j_1<\dots<j_r<j_{r+1}=k+1}\\
&~~~\left(\prod_{s=1}^rc_{a_s}^{(j_s)}\right)\prod_{s=0}^r\Ia(a_s;w_{j_s+1}\cdots w_{j_{s+1}-1};a_{s+1})\otimes\Im_\gamma(p;a_1,\dots,a_r;q),
\end{align*}
and the claim follows.
\epf

\subsection{Base Extension of $\cH_4$}
To discuss the motivic analogue of $F_2$-MZVs, the space $\cH_4$ is insufficient, and we have to consider the base extension.
\defi[Base Extension]
For a positive integer $N\in\bZge{1}$, put
\begin{align*}
\widetilde{\cH}_N:=\cH_N\otimes_\bQ \bQ_N.
\end{align*}
We also extend the following notations to $\widetilde\cH_N$: for $u\otimes c, u_1\otimes c_1, u_2\otimes c_2\in\widetilde\cH_N$,
\begin{itemize}
  \item coaction: $\Delta(u\otimes c):=\Delta(u)\otimes c$;
  \item product: $(u_1\otimes c_1)\cdot(u_2\otimes c_2):=u_1u_2\otimes c_1c_2$;
  \item period map: $\per(u\otimes c):=\per(u)\cdot c$.
\end{itemize}
In addition, we also consider the base extension of $\cU_N$ as  
\begin{align*}
\widetilde{\cU}_N:=\cU_N\otimes_\bQ \bQ_N.
\end{align*}
\edefi

By the conjecture that $\per: \cH_4\rightarrow\bC$ is injective (see Remark \ref{rem: period_conjecture}), it is expected that MMZVs of level $4$ lose no information of MZVs, and the following Proposition shows that if the conjecture is true, the base extension also loses no information.
\propo \label{prop: injectivity}
If the period map $\per:\cH_4\rightarrow\bC$ is injective, then its extension $\per:\widetilde\cH_4\rightarrow\bC$ is also injective.
\epropo
\pf
The element of $\widetilde\cH_4$ can be written as $u\otimes1+v\otimes\sqrt{-1}$ for some $u,v\in\cH_4$. If we assume $$\per(u\otimes1+v\otimes\sqrt{-1})=0,$$ then by multiplying $\per(u\otimes1-v\otimes\sqrt{-1})$, we have $$\per((u^2+v^2)\otimes1)=0.$$ Now, by the 
assumption, we have $u^2+v^2=0$. We will show that $u=v=0$ is only the case such that the elements $u,v\in\cH_4$ satisfy $u^2+v^2=0$. Here, we fix a non-canonical isomorphism $\Phi:\cH_4\rightarrow\cU_4$ in Theorem \ref{theo: str}, and we show that $u=v=0$ is only the case such that the elements $u,v\in\cH_4$ satisfy $\Phi(u)^{\shu2}+\Phi(v)^{\shu2}=0$ by contradiction. We assume $u\ne0$ or $v\ne0$. If $u=0$, then $\Phi(v)^{\shu2}=0$, and we have $v=0$ since $\cU_4$ is free. So it is sufficient if we consider the case where $u$ and $v$ together are non-zero. In addition, let $k_1$ and $k_2$ be the maximum degrees of $\Phi(u)$ and $\Phi(v)$ respectively, then the maximum degrees of $\Phi(u)^{\shu2}$ and $\Phi(v)^{\shu2}$ become $2k_1, 2k_2$. Therefore, we have $k_1=k_2$ since $\cU_4$ is free. Furthermore, if $u=\sum_{m=0}^ku^{(m)}$, $v=\sum_{m=0}^kv^{(m)}$~($u^{(m)}, v^{(m)}\in\cH_4^{(m)}$), then
\begin{align*}
0&=\Phi(u)\shu\Phi(u)+\Phi(v)\shu\Phi(v)\\
&=\Phi(u^{(k)})\shu\Phi(u^{(k)})+(\text{terms of degree less than } 2k-1),
\end{align*}
so it is sufficient if we consider the case of $u=u^{(k)}$ and $v=v^{(k)}\in\cH_4^{(k)}$. Let $$B_k=\left\{f_{j_1}\dots f_{j_r}\tau^l~\middle|~\sum_{s=1}^r j_s+l=k\right\}$$ be a basis of $\cU_4$, and define an order $\preceq$ on $B_k$ by
\begin{align*}
f_{j_1}\dots f_{j_r}\tau^l\preceq f_{j'_1}\dots f_{j'_s}\tau^{l'}\Longleftrightarrow &\text{there exists }m\ge1\text{ such that }\\
&j_1=j'_1,\dots j_{m-1}=j'_{m-1}, j_m\le j'_m.
\end{align*}
For $u'=u,v$ and $w\in B_k$, let $c_{u'}(w)\in\bQ$ be satisfying $\Phi(u')=\sum_{w\in B_k}c_{u'}(w)w$. Now, we assume $u,v\ne0$, and put $w_0=f_{j_1}\dots f_{j_r}\tau^l$ be the maximum element $w\in B_k$ such that $c_u(w)\ne0$ or $c_v(w)\ne0$. Then,
\begin{align*}
0&=\Phi(u^2+v^2)\\
&=\Phi(u)\shu\Phi(u)+\Phi(v)\shu\Phi(v)\\
&=(c_u(w_0)^2+c_v(w_0)^2)(f_{j_1}f_{j_1}\dots f_{j_r}f_{j_r}\tau^{2l})+(\text{terms less than } w_0).
\end{align*}
Since $c_u(w_0),c_v(w_0)\in\bQ$, we have $c_u(w_0)=c_v(w_0)=0$, and this contradicts the definition of $w_0$. Therefore, we have $u=v=0$, and the claim follows.
\epf

Now, we can define $F_2$-MMZVs as elements of $\widetilde\cH_4$. 

\defi \label{def: motivic-myzeta}
For $k_j\in\bZge{1}$ and $\varphi_j\in\{\omega_1,\omega\}~(j=1,\dots,d)$, we define MMZVs with respect to $F_2$ ($F_2$-MMZVs for short) by
\begin{align*}
\myzetam{k_1,\dots,k_d}{\varphi_1,\dots,\varphi_d}:=\Im_\dch(0;\theta_1\eta_0^{k_1-1}\dots\theta_d\eta_0^{k_d-1};1)\otimes(-\sqrt{-1})^l,
\end{align*}
where $l=\len_\omega(\bp)$ and 
$\theta_j=\begin{cases}
\eta, & \varphi_j=\omega,\\
\eta_1, & \varphi_j=\omega_1,
\end{cases}$ (see Definition \ref{def: phis} for the definition of $\omega,$ $\omega_1$, and Proposition \ref{prop: myzeta-level} for the definition of $\eta,$ $\eta_1$). Let $\MZV^{\mathfrak{m},(k)}_{F_2}$ be the $\bQ$-linear space spanned by $F_2$-MMZVs of weight $k$. Since Proposition \ref{prop: myzeta-level}, we have $$\per\left(\myzetam{k_1,\dots,k_d}{\varphi_1,\dots,\varphi_d}\right)=\myzeta{k_1,\dots,k_d}{\varphi_1,\dots,\varphi_d}.$$
\edefi

\subsection{Galois Action}
The Galois group $\Gal(\bQ_4/\bQ)$ acts on the category of mixed Tate motives $\mathrm{MTM}_{\mathcal{O}_4[1/4]}$ on $\mathcal{O}_4[1/4]$ induced by the automorphism as a ring of $\mathcal{O}_4[1/4]$ induced by the Galois action. In addition, the Galois action is induced on $\cH_4$ since $\cH_4$ is the coordinate ring of Tannakian fundamental group of $\mathrm{MTM}_{\mathcal{O}_4[1/4]}$. By Theorem \ref{theo: BD}, any elements of $\cH_4$ can be written by motivic iterated integrals, so in this paper, we define the Galois action explicitly for motivic iterated integrals.

\propo \label{prop: galois}
The Galois group $\Gal(\bQ_4/\bQ)=\{1,\sigma\}$ acts on $\widetilde\cH_4$ and this action is commutative with the production and the coaction $\Delta$. Moreover, $\sigma$ acts on the element written by motivic iterated integrals as
\begin{align*}
\sigma(\Im_\gamma(p;a_1,\dots,a_k;q)\otimes\epsilon)=\Im_{\sigma\gamma}({}^\sigma{p};{}^\sigma{a_1},\dots,{}^\sigma{a_k};{}^\sigma{q})\otimes{}^\sigma{\epsilon},
\end{align*}
where the path $\sigma\gamma: [0,1]\rightarrow \bC$ is defined by $\sigma\gamma(t):=\sigma(\gamma(t))$. In particular, if $\gamma$ is $\sigma$-invariant, we have
\begin{align*}
\sigma(\Im_\gamma(p;a_1,\dots,a_k;q)\otimes\epsilon)=\Im_{\gamma}(p;{}^\sigma{a_1},\dots,{}^\sigma{a_k};q)\otimes{}^\sigma{\epsilon}.
\end{align*}
\epropo

\rem
By Proposition \ref{prop: galois}, for $u\in\widetilde\cH_4$, we have $\per({}^\sigma u)={}^\sigma{\per(u)}$. In particular, it is equivalent that $u\in\widetilde\cH_4$ is $\sigma$-invariant and $\per(u)\in\bR$.
\erem
We also consider the Galois action on $\widetilde\cU_4$ induced by the non-canonical isomorphism $\Phi: \widetilde{\cH}_4\rightarrow\widetilde{\cU}_4$.

\thm \label{theo: 4to2}
Let $\left(\cD_1\cH_4^{(k)}\right)^\sigma$ be the $\sigma$-invariant part of $\cD_1\cH_4^{(k)}$, then we have
\begin{align*}
\left(\cD_1\cH_4^{(k)}\right)^\sigma=\cD_1\cH_2^{(k)}.
\end{align*}
\ethm
\pf
Let ${\rm MTM}_{\mathcal{O}_N[1/N]}$ be the category of mixed Tate motives over
$\mathcal{O}_N[1/N]$ (cf. {\cite[p.336]{Gl}}).
By using the tensor functor with
$h^0({\rm Spec}\hspace{0.5mm}\mathcal{O}_4[1/4])$, we have
a natural fully faithful functor from ${\rm MTM}_{\mathcal{O}_2[1/2]}$ to
${\rm MTM}_{\mathcal{O}_4[1/4]}$ (cf. \cite[2.1.6]{DG}).
It yields a natural injection $\iota:\mathcal H_2\hookrightarrow \mathcal H_4$ which preserves the depth grading and the motivic weight.
Therefore, the map $\iota$ also yields a natural injection
$\mathcal D_1 \mathcal H^{(k)}_2 \hookrightarrow \mathcal D_1 \mathcal H^{(k)}_4$ which commutes
with the inclusion $\mathcal D_1 \mathcal H^{(k)}_N\subset \mathcal H^{(k)}_N$ for
$N=2,4$.
Since any object in ${\rm MTM}_{\mathcal{O}_4[1/4]}$ is unramified outside 2 by
\cite[p.10,Proposition 1.8]{DG}, the claim follows from the Galois descend (see \cite[p.17, (2.16.2)]{DG}).
\epf

\section{Main Results and Their Proofs} \label{ss: proof}
\subsection{Preparations}
To state the claim of the main results, we ``twist'' the non-canonical isomorphism $\Phi: \cH_4\rightarrow\cU_4$ in Theorem \ref{theo: str}. Then, we can write the Galois action on the Hopf structure $\cU_4$ concisely.
\propo \label{prop: regPhi}
We can normalize the isomorphism $\widetilde\Phi: \cH_4\rightarrow\cU_4$ as
\begin{align*}
\widetilde\Phi(g_{2m-1})&=f_{2m-1}\\
\widetilde\Phi(h_{2m})&=f_{2m}
\end{align*}
for positive integers $m\ge1$, where 
\begin{align}
\begin{split}
g_k:=\frac{1}{2}\left(\zetamN{4}{k}{\sqrt{-1}}+\zetamN{4}{k}{-\sqrt{-1}}\right),\\
h_k:=\frac{1}{2}\left(\zetamN{4}{k}{\sqrt{-1}}-\zetamN{4}{k}{-\sqrt{-1}}\right). \label{eq: reg}
\end{split}
\end{align}
We also have $g_{2m}\in\bQ\tau^{2m}$ and $h_{2m-1}\in\bQ\tau^{2m-1}$. In addition, $\widetilde\Phi$ induces an isomorphism $\widetilde\Phi|_{\cH_2}:\cH_2\rightarrow\cU_2\subseteq\cU_4$.
\epropo
\pf
We fix an isomorphism $\Phi:\cH_4\rightarrow\cU_4$ in Theorem \ref{theo: str}, and normalize it to construct $\widetilde\Phi$. First, by Theorem \ref{theo: coaction}, Theorem \ref{theo: MII} (3), and
\begin{align*}
\zetamN{4}{k}{\pm\sqrt{-1}}=\Im_\dch(0;\mp\sqrt{-1},\{0\}^{k-1};1),
\end{align*}
we have $$\Delta\left(\zetamN{4}{k}{\pm\sqrt{-1}}\right)=1\otimes\zetamN{4}{k}{\pm\sqrt{-1}}+\rho\left(\zetamN{4}{k}{\pm\sqrt{-1}}\right)\otimes1.$$ Then, we have $\zetamN{4}{k}{-\sqrt{-1}}\in\Ker\widetilde\Delta$, and there exists $a_k, b_k\in\bQ$ such that
\begin{align*}
\zetamN{4}{k}{-\sqrt{-1}}=a_k\tau^k+b_k\zetamN{4}{k}{\sqrt{-1}}
\end{align*}
by Proposition \ref{prop: kernelM}. Furthermore, by Proposition \ref{prop: conjugate}, we have $b_k=(-1)^{k-1}$ and
\begin{align*}
g_{2m-1}&=\frac{a_{2m-1}}{2}\tau^{2m-1}+\zetamN{4}{2m-1}{\sqrt{-1}},\\
h_{2m}&=-\frac{a_{2m}}{2}\tau^{2m}+\zetamN{4}{2m}{\sqrt{-1}}
\end{align*}
for $m\ge1$. Therefore, $\Phi$ can be normalized as (\ref{eq: reg}) since $\zetamN{4}{k}{\sqrt{-1}}$ and $\tau^k$ are $\bQ$-linear independent. In addition, we have
\begin{align*}
g_{2m}&=\frac{a_{2m}}{2}\tau^{2m}\in\bQ\tau^{2m},\\
h_{2m-1}&=-\frac{a_{2m-1}}{2}\tau^{2m-1}\in\bQ\tau^{2m-1}.
\end{align*}
\epf

In the following, we fix the isomorphism by $\widetilde\Phi: \cH_4\rightarrow \cU_4$ given in Proposition \ref{prop: regPhi}.

\lem \label{lemm: sigma}
For a positive integer $k\in\bZge{1}$, put 
\begin{align*}
\left(\cU_4^{(k)}\right)^\sigma&:=\{u\in\cU_4^{(k)}~|~{}^\sigma u=u\},\\
\left(\cU_4^{(k)}\right)^{\sigma,-}&:=\{u\in\cU_4^{(k)}~|~{}^\sigma u=-u\}.
\end{align*}
Then, we have
\begin{align*}
\Ker\widetilde\Delta\cap\left(\cU_4^{(k)}\right)^\sigma&=\begin{cases}
\bQ f_k, & k \text{: odd},\\
\bQ \tau^k, & k \text{: even},
\end{cases}\\
\Ker\widetilde\Delta\cap\left(\cU_4^{(k)}\right)^{\sigma,-}&=\begin{cases}
\bQ \tau^k, & k \text{: odd},\\
\bQ f_k, & k \text{: even}.
\end{cases}
\end{align*}

\elem
\pf
By Lemma \ref{lemm: kernelH} and Proposition \ref{prop: regPhi}, we have $\Ker\widetilde\Delta~\cap~\cU_4^{(k)}=\bQ \widetilde\Phi(g_k)\oplus\bQ\widetilde\Phi(h_k)$, so we investigate the $\sigma$ action for $g_k, h_k$. By the definition of $g_k$ and Proposition \ref{prop: galois}, we have
\begin{align*}
{}^\sigma g_k&=\frac{1}{2}\left({}^\sigma\zetamN{4}{k}{\sqrt{-1}}+{}^\sigma\zetamN{4}{k}{-\sqrt{-1}}\right)\\
&=-\frac{1}{2}\left({}^\sigma\Im_\dch(0;-\sqrt{-1},\{0\}^{k-1};1)+{}^\sigma\Im_\dch(0;\sqrt{-1},\{0\}^{k-1};1)\right)\\
&=-\frac{1}{2}\left(\Im_\dch(0;\sqrt{-1},\{0\}^{k-1};1)+\Im_\dch(0;-\sqrt{-1},\{0\}^{k-1};1)\right)\\
&=\frac{1}{2}\left(\zetamN{4}{k}{\sqrt{-1}}+\zetamN{4}{k}{-\sqrt{-1}}\right)=g_k,
\end{align*}
then $g_k\in\left(\cU_4^{(k)}\right)^\sigma$. On the other hands, by the definition of $h_k$ and Proposition \ref{prop: galois}, we have
\begin{align*}
{}^\sigma h_k&=\frac{1}{2}\left({}^\sigma\zetamN{4}{k}{\sqrt{-1}}-{}^\sigma\zetamN{4}{k}{-\sqrt{-1}}\right)\\
&=-\frac{1}{2}\left({}^\sigma\Im_\dch(0;-\sqrt{-1},\{0\}^{k-1};1)-{}^\sigma\Im_\dch(0;\sqrt{-1},\{0\}^{k-1};1)\right)\\
&=-\frac{1}{2}\left(\Im_\dch(0;\sqrt{-1},\{0\}^{k-1};1)-\Im_\dch(0;-\sqrt{-1},\{0\}^{k-1};1)\right)\\
&=-\frac{1}{2}\left(\zetamN{4}{k}{\sqrt{-1}}-\zetamN{4}{k}{-\sqrt{-1}}\right)=-h_k,
\end{align*}
then $h_k\in\left(\cU_4^{(k)}\right)^{\sigma,-}$. Therefore, we have
\begin{align*}
\Ker\widetilde\Delta\cap\left(\cU_4^{(k)}\right)^\sigma&=\bQ \widetilde\Phi(g_k),\\
\Ker\widetilde\Delta\cap\left(\cU_4^{(k)}\right)^{\sigma,-}&=\bQ \widetilde\Phi(h_k).
\end{align*}
Also, Proposition \ref{prop: regPhi} yields
\begin{align*}
\widetilde\Phi(g_k)&=\begin{cases}
f_k, & k \text{: odd},\\
\tau^k, & k \text{: even},
\end{cases}\\
\widetilde\Phi(h_k)&=\begin{cases}
\tau^k, & k \text{: odd},\\
f_k, & k \text{: even},
\end{cases}
\end{align*}
then the claim follows.
\epf

Here, we prove that $g_k$ can be written explicitly by $\logm2$ or the motivic Riemann zeta values.

\propo \label{prop: gk}
For each positive integer $k\in\bZge{1}$, we have
\renewcommand{\arraystretch}{1.5} 
\begin{align*}
g_k=\begin{cases}
-\frac{1}{2}\log^\mathfrak{m}2, & k=1,\\
2^{-2k+1}(1-2^{k-1})\zetam(k), & k>1.
\end{cases}
\end{align*}
\renewcommand{\arraystretch}{1}
\epropo
\pf
By Theorem \ref{theo: 4to2}, we have $g_k\in\left(\cD_1\cH_4^{(k)}\right)^\sigma=\cD_1\cH_2^{(k)}$. Moreover, since $g_k\in\Ker\widetilde\Delta$ and by Proposition \ref{prop: kernelM}, there exists $c_k\in\bQ$ such that
\begin{align*}
g_k:=\frac{1}{2}\left(\zetamN{4}{k}{\sqrt{-1}}+\zetamN{4}{k}{-\sqrt{-1}}\right)=\begin{cases}
c_1\logm2, & k=1,\\
c_k\zetam(k), & k\ge2.
\end{cases}
\end{align*}
Now, by the direct calculation of $\per(g_k)$,
\renewcommand{\arraystretch}{1.5}
\begin{align*}
\per(g_k)&=-\sum_{n=1}^\infty\frac{(-1)^n}{(2n)^k}=\begin{cases}
-\frac{1}{2}\log2, & k=1,\\
-2^{-2k+1}(2^{k-1}-1)\zeta(k), & k>1,
\end{cases}
\end{align*}
\renewcommand{\arraystretch}{1}
and we have $c_1=-\frac{1}{2}$, $c_k=2^{-2k+1}(1-2^{k-1})~(k\ge2)$.
\epf

\subsection{Proof of Theorem \ref{theo: main1}}
Here, we give a proof of Theorem \ref{theo: main1}.
 
\lem \label{lemm: 1}
For positive integers $j_1,\dots,j_r\ge1$ and a non-negative integer $l\ge0$, the following equation holds:
\begin{align*}
{}^\sigma(f_{j_1}\dots f_{j_r}\tau^l)=(-1)^l\cdot{}^\sigma f_{j_1}\dots {}^\sigma f_{j_r}\tau^l.
\end{align*}
\elem
\pf
We prove the claim by induction on $r$. When $r=0$, since ${}^\sigma\tau=-\tau$, we have
\begin{align*}
{}^\sigma(\tau^l)=({}^\sigma\tau)^l=(-\tau)^l=(-1)^l\tau^l.
\end{align*}
When $r>0$, by the induction hypothesis, the equation (\ref{eq: coaU}), and Proposition \ref{prop: galois}, we have
\begin{align*}
\Delta'{}^\sigma(f_{j_1}\dots f_{j_r}\tau^l)&={}^\sigma(\Delta'(f_{j_1}\dots f_{j_r}\tau^l))\\
&={}^\sigma\left(\sum_{s=1}^rf_{j_1}\dots f_{j_s}\otimes f_{j_{s+1}}\dots f_{j_r}\tau^l)\right)\\
&=(-1)^l\sum_{s=1}^r{}^\sigma f_{j_1}\dots {}^\sigma f_{j_s}\otimes {}^\sigma f_{j_{s+1}}\dots {}^\sigma f_{j_r}\tau^l\\
&=(-1)^l\Delta'({}^\sigma f_{j_1}\dots {}^\sigma f_{j_r}\tau^l).
\end{align*}
Then, if $k=j_1\dots+j_r+l$, we have
\begin{align*}
{}^\sigma(f_{j_1}\dots f_{j_r}\tau^l)-(-1)^l\cdot{}^\sigma f_{j_1}\dots {}^\sigma f_{j_r}\cdot({}^\sigma\tau^l)\in\bQ\tau^k
\end{align*}
by Lemma \ref{lemm: kernelH}. Therefore, since $r>0$, we have
\begin{align*}
{}^\sigma(f_{j_1}\dots f_{j_r}\tau^l)=(-1)^l\cdot{}^\sigma f_{j_1}\dots {}^\sigma f_{j_r}\cdot({}^\sigma\tau^l).
\end{align*}
\epf

\pf[of Theorem \ref{theo: main1}]
First, for $k\in\bZge{1}$,
\begin{align*}
{}^\sigma f_{2k-1}&=\widetilde\Phi\left(\frac{1}{2}\left({}^\sigma\zetamN{4}{2k-1}{\sqrt{-1}}+{}^\sigma\zetamN{4}{2k-1}{-\sqrt{-1}}\right)\right)\\
&=\widetilde\Phi\left(\frac{1}{2}\left(\zetamN{4}{2k-1}{-\sqrt{-1}}+\zetamN{4}{2k-1}{\sqrt{-1}}\right)\right)=f_{2k-1},\\
{}^\sigma f_{2k}&=\widetilde\Phi\left(\frac{1}{2}\left({}^\sigma\zetamN{4}{2k}{\sqrt{-1}}-{}^\sigma\zetamN{4}{2k}{-\sqrt{-1}}\right)\right)\\
&=\widetilde\Phi\left(\frac{1}{2}\left(\zetamN{4}{2k}{-\sqrt{-1}}-\zetamN{4}{2k}{\sqrt{-1}}\right)\right)=-f_{2k},
\end{align*}
then we have ${}^\sigma f_j=(-1)^{j-1}f_j$ for any $j\in\bZge{1}$. Therefore, for $f_{j_1}\dots f_{j_r}\tau^l\otimes\epsilon~~(j_1+\dots+j_r+l=k,~\epsilon\in\{1, \sqrt{-1}\})$, we have
\begin{align*}
{}^\sigma(f_{j_1}\dots f_{j_r}\tau^l\otimes\epsilon)&=(-1)^l\cdot{}^\sigma f_{j_1}\dots {}^\sigma f_{j_r}\tau^l\otimes{}^\sigma\epsilon\\
&=(-1)^{l+(j_1-1)+\dots+(j_r-1)}f_{j_1}\dots f_{j_r}\tau^l\otimes{}^\sigma\epsilon\\
&=(-1)^{k-r}f_{j_1}\dots f_{j_r}\tau^l\otimes{}^\sigma\epsilon
\end{align*}
by Lemma \ref{lemm: 1}. Then, we have
\begin{align*}
f_{j_1}\dots f_{j_r}\tau^l\otimes\epsilon\in\left(\widetilde\cH^{(k)}\right)^\sigma\Longleftrightarrow \epsilon=\begin{cases}
1, & k-r \text{: even},\\
\sqrt{-1}, & k-r \text{: odd}.
\end{cases}
\end{align*}
\epf

\subsection{Corollaries of Theorem \ref{theo: main1}} \label{ss: cors}

By Theorem \ref{theo: main1}, we can identify the dimension of $\left(\cH_4^{(k)}\right)^\sigma$ and the upper bound of the dimension of $\MZV_{F_2}^{(k)}$.
\cor \label{cor: dimension}
For each non-negative integer $k\ge0$, the following equation holds:
\begin{align*}
\dim_\bQ \left(\widetilde\cH_4^{(k)}\right)^\sigma=2^k.
\end{align*}
In particular, we have 
\begin{align*}
\dim_\bQ \MZV_{F_2}^{(k)}\le2^k.
\end{align*}
\ecor
\pf
Put $d_k:=\dim_\bQ\left(\widetilde\cH_4^{(k)}\right)^\sigma$, then
\begin{align*}
\sum_{k\ge0}d_kt^k
&=\left\{\sum_{r=0}^\infty(t+t^2+\dots)^r\right\}(1+t+t^2+\dots)\\
&=\frac{1}{1-(t+t^2+t^3+\dots)}\cdot\frac{1}{1-t}\\
&=\frac{1}{1-t\cdot(1-t)^{-1}}\cdot\frac{1}{1-t}=\frac{1}{1-2t}.
\end{align*}
Therefore, we have $d_k=2^k$. Furthermore, we have $\dim_\bQ\MZV_{F_2}^{(k)}\le2^k$ since $\MZV_{F_2}^{(k)}\subseteq\per\left(\left(\widetilde\cH_4^{(k)}\right)^\sigma\right)$.
\epf

We can also show that $\MZV^{\mathfrak{m},(k)}_{F_2}$ has a direct sum decomposition with respect to the parity of $\len_\omega(\bp)$.
\cor \label{cor: 2}
For each positive integer $k\ge1$, the following decomposition holds:
\begin{align*}
\MZV^{\mathfrak{m},(k)}_{F_2}&=\MZV^{\mathfrak{m},(k)}_{F_2,\mathrm{odd}}\oplus\MZV^{\mathfrak{m},(k)}_{F_2,\mathrm{even}},
\end{align*}
where
\begin{align*}
\MZV^{\mathfrak{m},(k)}_{F_2,\mathrm{odd}}&:=\left\langle\myzetam{k_1,\dots,k_d}{\varphi_1,\dots,\varphi_d}~\middle|~\substack{k_1+\dots+k_d=k,\\ \len_\omega(\bp)\mathrm{: odd}}\right\rangle_\bQ,\\
\MZV^{\mathfrak{m},(k)}_{F_2,\mathrm{even}}&:=\left\langle\myzetam{k_1,\dots,k_d}{\varphi_1,\dots,\varphi_d}~\middle|~\substack{k_1+\dots+k_d=k,\\ \len_\omega(\bp)\mathrm{: even}}\right\rangle_\bQ.
\end{align*}
Also, we have
\begin{align*}
\dim_\bQ\MZV^{\mathfrak{m},(k)}_{F_2,\mathrm{odd}}&\le2^{k-1},\\
\dim_\bQ\MZV^{\mathfrak{m},(k)}_{F_2,\mathrm{even}}&\le2^{k-1}.
\end{align*}
\ecor

\rem
By Proposition \ref{prop: injectivity}, we have
\begin{align*}
\dim_\bQ\MZV^{(k)}_{F_2,\mathrm{odd}}&\le 2^{k-1},\\
\dim_\bQ\MZV^{(k)}_{F_2,\mathrm{even}}&\le 2^{k-1}.
\end{align*}
Also, if $\per: \cH_4\rightarrow\bC$ is injective (see Remark \ref{rem: period_conjecture}), $\MZV^{(k)}_{F_2}$ also has the same decomposition:
\begin{align*}
\MZV^{(k)}_{F_2}&=\MZV^{(k)}_{F_2,\mathrm{odd}}\oplus\MZV^{(k)}_{F_2,\mathrm{even}},
\end{align*}
where
\begin{align*}
\MZV^{(k)}_{F_2,\mathrm{odd}}&:=\left\langle\myzeta{k_1,\dots,k_d}{\varphi_1,\dots,\varphi_d}~\middle|~\substack{k_1+\dots+k_d=k,\\ \len_\omega(\bp)\mathrm{: odd}}\right\rangle_\bQ,\\
\MZV^{(k)}_{F_2,\mathrm{even}}&:=\left\langle\myzeta{k_1,\dots,k_d}{\varphi_1,\dots,\varphi_d}~\middle|~\substack{k_1+\dots+k_d=k,\\ \len_\omega(\bp)\mathrm{: even}}\right\rangle_\bQ.
\end{align*}
\erem

\pf[of Corollary \ref{cor: 2}]
First, for $\widetilde\Phi\left(\widetilde\cH_4^{(k)}\right)^\sigma$, we have the following decomposition:
\begin{align*}
\widetilde\Phi\left(\widetilde\cH_4^{(k)}\right)^\sigma
&=\left\langle f_{j_1}\dots f_{j_r}\tau^l\otimes1~|~j_1+\dots+j_r+l=k, k-r\text{: even}\right\rangle_\bQ\\
&\oplus\left\langle f_{j_1}\dots f_{j_r}\tau^l\otimes\sqrt{-1}~|~j_1+\dots+j_r+l=k, k-r\text{: odd}\right\rangle_\bQ.
\end{align*}
Also, we have the following equations:
\begin{align*}
\left\langle f_{j_1}\dots f_{j_r}\tau^l\otimes1~|~j_1+\dots+j_r+l=k, k-r\text{: even}\right\rangle_\bQ&=\widetilde\Phi\left(\left(\cH_4^{(k)}\right)^\sigma\right),\\
\left\langle f_{j_1}\dots f_{j_r}\tau^l\otimes\sqrt{-1}~|~j_1+\dots+j_r+l=k, k-r\text{: odd}\right\rangle_\bQ&=\widetilde\Phi\left(\left(\cH_4^{(k)}\otimes\sqrt{-1}\bQ\right)^\sigma\right).
\end{align*}
On the other hands, by definition of $\myzetam{k_1,\dots,k_d}{\varphi_1,\dots,\varphi_d}\in\MZV_{F_2}^{\mathfrak{m},(k)}$, we have
\begin{align*}
\len_\omega(\boldsymbol{\varphi}) \text{: even} &\Longrightarrow \myzetam{k_1,\dots,k_d}{\varphi_1,\dots,\varphi_d}\in\left(\cH_4^{(k)}\right)^\sigma,\\
\len_\omega(\boldsymbol{\varphi}) \text{: odd} &\Longrightarrow \myzetam{k_1,\dots,k_d}{\varphi_1,\dots,\varphi_d}\in\left(\cH_4^{(k)}\otimes\sqrt{-1}\bQ\right)^\sigma,
\end{align*}
then
\begin{align*}
\MZV^{\mathfrak{m},(k)}_{F_2,\mathrm{even}}&\subseteq\widetilde\Phi\left(\left(\cH_4^{(k)}\right)^\sigma\right),\\
\MZV^{\mathfrak{m},(k)}_{F_2,\mathrm{odd}}&\subseteq\widetilde\Phi\left(\left(\cH_4^{(k)}\otimes\sqrt{-1}\bQ\right)^\sigma\right),
\end{align*}
and the claim follows.
\epf

\subsection{Preparation for the Proof of Theorem \ref{theo: main2}}
Now, we show some lemmas for the proof of Theorem \ref{theo: main2}. Recall that $\myzetam{k}{\omega}$ is defined by
\begin{align*}
\myzetam{k}{\omega}=\Im_\dch(0;\eta\eta_0^{k-1};1)\otimes(-\sqrt{-1})\in\widetilde\cH_4,
\end{align*}
where $\eta=e_i-e_{-i}$ and $\eta_0=e_0-e_i-e_{-i}\in W_{\widetilde\mu_4}^{(1)}$. Put $$u_k=\Im_\dch(0;\eta\eta_0^{k-1};1)\in\cH_4.$$ We prove Theorem \ref{theo: main2} by showing that the equation holds in $\widetilde\cU_4$ if we send the both side of the equation (\ref{eq: main2}) by the isomorphism $\widetilde\Phi: \widetilde\cH_4\rightarrow\widetilde\cU_4$. The following three lemmas (Lemma \ref{lemm: 3}-\ref{lemm: 4}) are the calculations appearing in the coaction of $u_k$, and in Lemma \ref{lem: 5}, we describe where the right hand side of the equation (\ref{eq: main2}) sent into $\widetilde\cU_4$ by $\widetilde\Phi$.

\lem \label{lemm: 3}
For each non-negative integer $r\in\bZge{1}$, we have
\begin{enumerate}
  \item $\Ia(-i;\eta_0^r;i)=0$,
  \item $\Ia(\pm i;\eta_0^r;1)=0$,
  \item $\Ia(0;\eta\eta_0^{r-1};i)=-\Ia(0;\eta\eta_0^{r-1};-i)$.
\end{enumerate}
\elem
\pf
First, we prove (1) in the case of $r=1$, i.e. $\Ia(-i;\eta_0;i)=0$. By applying Theorem \ref{theo: MII} (7) for $c=-1$, we have
\begin{align*}
\Ia(-i;\eta_0;i)=\Ia(-i;e_0-e_i-e_{-i};i)=\Ia(i;e_0-e_{-i}-e_i;-i)=\Ia(i;\eta_0;-i).
\end{align*}
On the other hands, by Theorem \ref{theo: MII} (5), we have
\begin{align*}
\Ia(-i;\eta_0;i)=-\Ia(i;e_0-e_i-e_{-i};-i)=-\Ia(i;\eta_0;-i).
\end{align*}
By comparing the two equations above, we have $\Ia(-i;\eta_0;i)=0$. If $r>1$, by Lemma \ref{lemm: exp}, we have
\begin{align*}
\Ia(-i;\eta_0^r;i)=\frac{1}{r!}\Ia(-i;\eta_0^{\shu r};i)=\frac{1}{r!}\Ia(-i;\eta_0;i)^r=0.
\end{align*}
Secondaly, we prove (2) in the case of $r=1$, i.e. $\Ia(\pm i;\eta_0;1)=0$. By applying Theorem \ref{theo: MII} (7) for $c=i$, we have
\begin{align*}
\Ia(\pm i;0;1)=\Ia(\mp 1;0;i)=\Ia(\mp i;0;-1)=\Ia(\pm1;0;-i).
\end{align*}
Then, by Proposition \ref{prop: pathconnect}, we have
\begin{align*}
\Ia(\pm i;0;1)&=\frac{1}{4}\left\{\Ia(\pm i;0;1)+\Ia(\mp 1;0;i)+\Ia(\mp i;0;-1)+\Ia(\pm1;0;-i)\right\}\\
&=\frac{1}{4}\Ia(\pm i;0;\pm i)=0.
\end{align*}
Now, we calcurate
\begin{align*}
\Ia(\pm i;i;1)+\Ia(\pm i;-i;1)=\Ia(\pm i;i;0)+\Ia(0;i;1)+\Ia(\pm i;-i;0)+\Ia(0;-i;1).
\end{align*}
By Theorem \ref{theo: MII} (5), (7), we have
\begin{align*}
\Ia(-i;i;0)&=\Ia(1;-1;0)=-\Ia(0;-1;1),\\
\Ia(i;-i;0)&=\Ia(1;-1;0)=-\Ia(0;-1;1).
\end{align*} 
Moreover, by Theorem \ref{theo: MII} (4), (5), (7), we have
\begin{align*}
\Ia(i;i;0)=\Ia(1;1;0)=-\Ia(0;1;1)&=0,\\
\Ia(-i;-i;0)=\Ia(1;1;0)=\Ia(0;1;1)&=0.
\end{align*} 
In addition, by Proposition \ref{prop: distribution},
\begin{align*}
\Ia(0;i;1)+\Ia(0;-i;1)=\Ia(0;-1;1),
\end{align*}
and we have $\Ia(\pm i;\eta_0;1)=0$. If $r>1$, by Lemma \ref{lemm: exp}
\begin{align*}
\Ia(\pm i;\eta_0^r;1)=\frac{1}{r!}\Ia(\pm i;\eta_0^{\shu r};1)=\frac{1}{r!}\Ia(\pm i;\eta_0;1)^r=0.
\end{align*}
Now, we prove (3). By applying Theorem \ref{theo: MII} (7) for $c=-1$, since $\eta=e_i-e_{-i}$ is changed to $-\eta$, we have
\begin{align*}
\Ia(0;\eta\eta_0^{r-1};i)=-\Ia(0;\eta\eta_0^{r-1};-i).
\end{align*}
\epf

The next is the key lemma for the proof of Theorem \ref{theo: main2}.
\newcommand{\I}{\mathrm{I}}
\lem \label{lemm: coaction_uv}
Let $k\ge1$ be a positive integer. Then, for the coaction $\Delta: \widetilde\cU_4\rightarrow\widetilde\cV_4\otimes_\bQ\widetilde\cU_4$, we have
\begin{align}
\label{eq: coaction_u} \widetilde\Delta(u_k)=-\sum_{r=1}^{k-1}v_{k-r}\otimes u_r,
\end{align}
where
\begin{align*}
v_k:=\Ia(0;\eta\eta_0^{k-1};i)\in\cA_4^{(k)}.
\end{align*}
Furthermore, we have
\begin{align}
\label{eq: coaction_v}\widetilde\Delta(v_k)=-\sum_{r=1}^{k-1}v_{k-r}\otimes v_r.
\end{align}
\elem

\pf
We prove the equations (\ref{eq: coaction_u}) and (\ref{eq: coaction_v}) at the same time. We put $u_k(q)=\begin{cases}
u_k, & q=1,\\
v_k, & q=i,
\end{cases}$ for $q\in\{1,i\}$, and $$\I(0;a_1,\dots,a_r;q):=\begin{cases}
\Im_\dch(0;a_1,\dots,a_r;1), & q=1,\\
\Ia(0;a_1,\dots,a_r;i), & q=i,
\end{cases}$$ for $a_1,\dots,a_r\in\widetilde\mu_4$. Then, by Proposition \ref{prop: coaction2}, we have
\begin{align} \label{eq: coaction_uv}
\widetilde\Delta u_k(q)&=\widetilde\Delta\Ia(0;\eta\eta_0^{k-1};q)\\
\notag &=\sum_{r=1}^{k-1}\sum_{(a_1,\dots,a_r)\in(\widetilde\mu_4)^r}\sum_{0=j_0<j_1<\dots<j_r<j_{r+1}=k+1}\\
\notag &\left(\prod_{s=1}^rc_{a_s}^{(j_s)}\right)\cdot\left(\prod_{s=0}^r\Ia(a_s;w_{j_s+1},\dots,w_{j_{s+1}-1};a_{s+1})\right)\otimes\I(0;a_1,\dots,a_r;q).
\end{align}
Here, since $w_j=\begin{cases}
\eta=e_i-e_{-i}, & j=1,\\
\eta_0=e_0-e_i-e_{-i}, & j>1,
\end{cases}$ the coefficients $c_a^{(j)}$ become
\begin{align*}
c_{a}^{(j)}&=\begin{cases}
0, & a=\pm1 \text{ and } j\ge1,\\
0, & a=0 \text{ and } j=1,\\
1, & a=0 \text{ and } j>1,\\
1, & a=i \text{ and } j=1,\\
-1, & a=i \text{ and } j>1,\\
-1, & a=-i \text{ and } j\ge1,
\end{cases}
\end{align*}
for each $a\in\widetilde\mu_4$ and $j=1,\dots,k$. Note that $c_a^{(j)}$ is independent of $j>1$. Now, we fix $r=1,\dots,k-1$ in the equation (\ref{eq: coaction_uv}), and consider two different cases respectively: (i) the terms of $j_1=1$ and (ii) the terms of $j_1>1$.

\noindent
(i) When $j_1=1$, since $c_{0}^{(1)}=c_{1}^{(1)}=c_{-1}^{(1)}=0$, the terms without the cases of $a_1=\pm i$ become $0$. Also, if $s>1$, we have $j_s>1$ and $c_{a_s}^{(j_s)}=c_{a_s}^{(2)}$. Then,
\begin{align*}
&\sum_{1<j_2<\dots<j_r<j_{r+1}=k+1}\left(\prod_{s=1}^rc_{a_s}^{(j_s)}\right)\cdot\left(\prod_{s=0}^r\Ia(a_s;w_{j_s+1},\dots,w_{j_{s+1}-1};a_{s+1})\right)\otimes\I(0;a_1,\dots,a_r;q)\\
&=\left(c_{a_1}^{(1)}\prod_{s=2}^rc_{a_s}^{(2)}\right)\sum_{1<j_2<\dots<j_r<j_{r+1}=k+1}\left(\prod_{s=1}^r\Ia(a_s;\eta_0^{j_{s+1}-j_s-1};a_{s+1})\right)\otimes\I(0;a_1,\dots,a_r;q)
\end{align*}
for each $a_1\in\{\pm\sqrt{-1}\}$ and $a_2,\dots,a_r\in\{0,\pm\sqrt{-1}\}$. Now, we put $j_s'=j_s-s$ for $s=1,\dots,r$, then the range of sum is $0=j_1'\le j_2' \le\dots\le j_r'\le j_{r+1}'=k-r$. Also, $\sum_{s=1}^r(j_{s+1}'-j_s')=k-r$ is independent of the choice of $j_2,\dots,j_r$, then we have
\begin{align*}
&\sum_{1=j_1<j_2<\dots<j_r<j_{r+1}=k+1}\left(\prod_{s=1}^r\Ia(a_s;\eta_0^{j_{s+1}-j_s-1};a_{s+1})\right)\\
&=\sum_{0=j_1'\le j_2' \le\dots\le j_r'\le j_{r+1}'=k-r}\left(\prod_{s=1}^r\Ia(a_s;\eta_0^{j_{s+1}'-j_s'};a_{s+1})\right)\\
&=\Ia(a_1;\eta_0^{k-r};q)
\end{align*}
by Proposition \ref{prop: pathconnect}. In addition, by $r<k$ and Lemma \ref{lemm: 3} (1) and (2), we have
\begin{align*}
\sum_{1=j_1<j_2<\dots<j_r<j_{r+1}=k+1}\left(\prod_{s=1}^r\Ia(a_s;\eta_0^{j_{s+1}-j_s-1};a_{s+1})\right)=0.
\end{align*}

\noindent
(ii) Fix $j_1$ as $1<j_1\le k-r+1$. Then, for $a_s\in\{0,\pm\sqrt{-1}\}$ ($s=1,\dots,r$), $c_{a_s}^{j_s}=\begin{cases}
1, & a_s=0,\\
-1, & a_s=\pm\sqrt{-1}
\end{cases}$ is independent of the choice of $j_s$. Therefore, by Proposition \ref{prop: pathconnect}, we have
\begin{align*}
&\sum_{j_1<j_2<\dots<j_r<j_{r+1}=k+1}\left(\prod_{s=1}^rc_{a_s}^{(j_s)}\right)\cdot\left(\prod_{s=0}^r\Ia(a_s;w_{j_s+1},\dots,w_{j_{s+1}-1};a_{s+1})\right)\otimes\I(0;a_1,\dots,a_r;q)\\
&=\left(\prod_{s=1}^rc_{a_s}^{(j_s)}\right)\cdot\Ia(0;\eta\eta_0^{j_1-2};a_1)\\
&~~~~\times\sum_{j_1<j_2<\dots<j_r<j_{r+1}=k+1}\left(\prod_{s=1}^r\Ia(a_s;\eta_0^{j_{s+1}-j_s-1};a_{s+1})\right)\otimes\I(0;a_1,\dots,a_r;q)
\end{align*}
for each $(a_1,\dots,a_r)\in\{0,\pm\sqrt{-1}\}^r$. Now, we put $j_s'=j_s-s-j_1+1$ for each $s=1,\dots,r$, then the range of sum is $0=j_1'\le j_2' \le\dots\le j_r'\le j_{r+1}'=k-r-j_1+1$. Also, $\sum_{s=1}^r(j_{s+1}'-j_s')=k-r-j_1+1$ is independent of the choice of $j_2,\dots,j_r$, then we have
\begin{align*}
&\sum_{j_1<j_2<\dots<j_r<j_{r+1}=k+1}\left(\prod_{s=1}^r\Ia(a_s;\eta_0^{j_{s+1}-j_s-1};a_{s+1})\right)\\
&=\sum_{0=j_1'\le j_2'\le\dots\le j_r\le j_{r+1}=k-r-j_1+1}\left(\prod_{s=1}^r\Ia(a_s;\eta_0^{j_{s+1}-j_s-1};a_{s+1})\right)\\
&=\Ia(a_1;\eta_0^{k-r-j_1+1};i)
\end{align*}
by Proposition \ref{prop: pathconnect}. Now, by Lemma \ref{lemm: 3}, we have $\Ia(a_1;\eta_0^{k-r-j_1+1};i)=\begin{cases}
0, & r<k-j_1+1,\\
1, & r=k-j_1+1.
\end{cases}$ In addition, for $a_1=0$, we have $\Ia(0;\eta\eta_0^{j_1-2};a_1)=0$. Then, the terms become $0$ without the case of $r=k-j_1+1$ and $a_1=\pm\sqrt{-1}$. Therefore, by calculating sum in the range of $a_1\in\{\pm\sqrt{-1}\},$ and $a_2,\dots,a_r\in\{0,\pm\sqrt{-1}\}$, we have
\begin{align*}
&\sum_{(a_1,\dots,a_r)\in(\widetilde\mu_4)^r}\sum_{j_1<j_2<\dots<j_r<j_{r+1}=k+1}\left(\prod_{s=1}^rc_{a_s}^{(j_s)}\right)\\
&~~~\times\left(\prod_{s=0}^r\Ia(a_s;w_{j_s+1},\dots,w_{j_{s+1}-1};a_{s+1})\right)\otimes\I(0;a_1,\dots,a_r;q)\\
&=\sum_{\substack{a_1\in\{\pm\sqrt{-1}\},\\ a_2,\dots,a_r\in\{0,\pm\sqrt{-1}\}}}\left(\prod_{s=1}^rc_{a_s}^{(j_s)}\right)\Ia(0;\eta\eta_0^{j_1-2};a_1)\otimes \I(0;e_{a_1}e_{a_2}\dots e_{a_r};q)\\
&=-\sum_{a_2,\dots,a_r\in\{0,\pm\sqrt{-1}\}}\left(\prod_{s=2}^rc_{a_s}^{(j_s)}\right)\Ia(0;\eta\eta_0^{k-r-1};i)\otimes \I(0;e_{i}e_{a_2}\dots e_{a_r};q)\\
&~~~-\sum_{a_2,\dots,a_r\in\{0,\pm\sqrt{-1}\}}\left(\prod_{s=2}^rc_{a_s}^{(j_s)}\right)\Ia(0;\eta\eta_0^{k-r-1};-i)\otimes \I(0;e_{-i}e_{a_2}\dots e_{a_r};q).
\end{align*}
Moreover, by Lemma \ref{lemm: 3} (2), $\Ia(0;\eta\eta_0^{k-r-1};i)=-\Ia(0;\eta\eta_0^{k-r-1};-i)$, then
\begin{align*}
&-\sum_{a_2,\dots,a_r\in\{0,\pm\sqrt{-1}\}}\left(\prod_{s=2}^rc_{a_s}^{(j_s)}\right)\Ia(0;\eta\eta_0^{k-r-1};i)\otimes \I(0;e_{i}e_{a_2}\dots e_{a_r};q)\\
&-\sum_{a_2,\dots,a_r\in\{0,\pm\sqrt{-1}\}}\left(\prod_{s=2}^rc_{a_s}^{(j_s)}\right)\Ia(0;\eta\eta_0^{k-r-1};-i)\otimes \I(0;e_{-i}e_{a_2}\dots e_{a_r};q)\\
&=-\sum_{a_2,\dots,a_r\in\{0,\pm\sqrt{-1}\}}\Ia(0;\eta\eta_0^{k-r-1};i)\otimes \left(\prod_{s=2}^rc_{a_s}^{(j_s)}\right)\cdot\Bigl\{\I(0;e_{i}e_{a_2}\dots e_{a_r};q)-\I(0;e_{-i}e_{a_2}\dots e_{a_r};q)\Bigr\}.
\end{align*}
Furthermore, since $c_i^{(1)}=1$ and $c_{-i}^{(1)}=-1$, we have
\begin{align*}
&\left(\prod_{s=2}^rc_{a_s}^{(j_s)}\right)\cdot\Bigl\{\I(0;e_{i}e_{a_2}\dots e_{a_r};q)-\I(0;e_{-i}e_{a_2}\dots e_{a_r};q)\Bigr\}\\
&=\sum_{\substack{a_1\in\{\pm\sqrt{-1}\},\\ a_2,\dots,a_r\in\{0,\pm\sqrt{-1}\}}}c_{a_1}^{(1)}\cdot\left(\prod_{s=2}^rc_{a_s}^{(j_s)}\right) \cdot \I(0;e_{a_1}e_{a_2}\dots e_{a_r};q)\\
&=\I(0;\eta\eta_0^{r-1};q)\\
&=u_r(q).
\end{align*}
Therefore,
\begin{align*}
&\sum_{(a_1,\dots,a_r)\in(\widetilde\mu_4)^r}\sum_{j_1<j_2<\dots<j_r<j_{r+1}=k+1}\left(\prod_{s=1}^rc_{a_s}^{(j_s)}\right)\\
&~~~\times\left(\prod_{s=0}^r\Ia(a_s;w_{j_s+1},\dots,w_{j_{s+1}-1};a_{s+1})\right)\otimes\I(0;a_1,\dots,a_r;q)\\
&=-v_{k-r}\otimes u_r(q).
\end{align*}

Now, by calculating the sum in the range of $r=1,\dots,k-1$ for each (i) and (ii), and the equation (\ref{eq: coaction_uv}), we have
\begin{align*}
\widetilde\Delta(u_k(q))=-\sum_{r=1}^{k-1}v_{k-r}\otimes u_r(q).
\end{align*}
\epf

\lem \label{lemm: 4}
For each positive odd number $j\ge1$, there exists $\beta'_j\in\bQ$ such that the following equation holds:
\begin{align}
\widetilde{\Phi}(v_k)=(-1)^{k-1}\sum_{\mathbf{l}\in \bI_k^{\od}}\underset{j\ge1 \text{: odd}}{\Shu}(\beta_j'f_j)^{l_j} \label{eq: lem4}
\end{align}
for any $k\in\bZge{1}$, where 
\begin{align*}
\bI_k^{\od}:=\left\{(l_1,l_2,\dots,l_s)\in\bI_k~|~l_j=0~\text{if }j\text{: even}\right\}.
\end{align*}
\elem

\pf
We prove the claim by induction on $k$. When $k=1$, by Theorem \ref{theo: MII} (7), we have $$v_1=\Ia(0;\eta;i)=\Ia(0;e_0-e_1-e_{-1};1)\in\cH_2^{(1)}.$$ In particular, by Theorem \ref{theo: main1}, \ref{theo: str}, and since ${}^\sigma v_1=v_1$, we have $\Phi(v_1)\in\bQ f_1$, and $\Phi(v_1)=\beta_1'f_1$ for some $\beta_1'\in\bQ$. When $k>1$, we have
\begin{align}
\widetilde\Delta(v_k)&=-\sum_{r=1}^{k-1}v_{k-r}\otimes v_r, \notag\\
\Delta(v_k)&=-\sum_{r=1}^{k-1}v_{k-r}\otimes v_r+1\otimes v_k+v_k\otimes1, \label{eq: 1}
\end{align}
by the equation (\ref{eq: coaction_v}) in Lemma \ref{lemm: coaction_uv}. On the other hands, for the right hand side of the equation (\ref{eq: lem4}), we put
\begin{align*}
\widetilde v'_{k}:=(-1)^{k-1}\sum_{\mathbf{l}\in \bI_k^{\od}}\underset{j\ge1 \text{: odd}}{\Shu}(\beta'_jf_j)^{l_j}\in\cV_2^{(k)},
\end{align*}
then by the explicit formula of coaction on $\cU_4$ (the equation (\ref{eq: coaU})) and the commutativity of $\Delta$ and $\shu$,
\begin{align*}
\Delta(\widetilde v'_{k})&=(-1)^{k-1}\sum_{\mathbf{l}\in \bI_k^{\od}}\underset{j\ge1 \text{: odd}}{\Shu}\Delta((\beta'_jf_j)^{l_j})\\
&=(-1)^{k-1}\sum_{\mathbf{l}\in \bI_k^{\od}}\underset{j\ge1 \text{: odd}}{\Shu}\sum_{m_j=0}^{l_j}(\beta'_jf_j)^{k_j-m_j}\otimes (\beta'_jf_j)^{m_j}\\
&=(-1)^{k-1}\sum_{\mathbf{l}\in \bI_k^{\od}}\sum_{m_1+n_1=l_1}\dots\sum_{m_{|\mathbf{l}|}+n_{|\mathbf{l}|}=l_{|\mathbf{l}|}}\left(\underset{j\ge1 \text{: odd}}{\Shu}((\beta'_jf_j)^{m_j}\right)\otimes \left(\underset{j\ge1 \text{: odd}}{\Shu}(\beta_jf_j)^{n_j}\right).
\end{align*}
Now, we put $m=m_1+3m_3+\dots+|\bl|m_{|\bl|}$ and $n=n_1+3n_3+\dots+|\bl|n_{|\bl|}=k-m$, and change the order of the sum: calculating the sum of the terms for each $m=0,\dots,k$. Then, we have
\begin{align*}
\Delta(\widetilde v_{k}')&=-\sum_{m+n=k}\left((-1)^{m-1}\sum_{\mathbf{m}\in \bI_m^{\od}}\underset{j\ge1 \text{: odd}}{\Shu}(\beta'_jf_j)^{m_j}\right)\otimes\left((-1)^{n-1}\sum_{\mathbf{n}\in \bI_{n}^{\od}}\underset{j\ge1 \text{: odd}}{\Shu}(\beta'_jf_j)^{n_j}\right).
\end{align*}
Now, by the induction hypothesis, we have
\begin{align*}
\widetilde\Phi(v_m)&=(-1)^{m-1}\sum_{\mathbf{m}\in \bI_m^{\od}}\underset{j\ge1 \text{: odd}}{\Shu}(\beta'_jf_j)^{m_j}\\
\widetilde\Phi(v_n)&=(-1)^{n-1}\sum_{\mathbf{n}\in \bI_n^{\od}}\underset{j\ge1 \text{: odd}}{\Shu}(\beta'_jf_j)^{n_j}
\end{align*}
for $0<m, n<k$, then
\begin{align} \label{eq: 2}
\Delta(\widetilde v_{k}')=-\sum_{\substack{m+n=k\\m,n>0}}\widetilde\Phi(v_m)\otimes\widetilde\Phi(v_{n})+1\otimes\widetilde v'_k+\widetilde v'_k\otimes1.
\end{align}

Here, $\widetilde\Phi$ is an isomorphism as Hopf algebras, and commutative with $\Delta$, then
\begin{align*}
\Delta(\widetilde\Phi(v_k)-\widetilde v_k')
&=(\widetilde\Phi\otimes\widetilde\Phi)(\Delta(v_k))-\Delta(\widetilde v'_k)\\
&=1\otimes\widetilde\Phi(v_k)+\widetilde\Phi(v_k)\otimes1-1\otimes\widetilde v'_k-\widetilde v'_k\otimes1\\
&=1\otimes(\widetilde\Phi(v_k)-\widetilde v'_k)+(\widetilde\Phi(v_k)-\widetilde v'_k)\otimes1
\end{align*}
by the equations (\ref{eq: 1}) and (\ref{eq: 2}). Therefore, we have $\widetilde\Phi(v_k)-\widetilde v_k'\in\Ker\widetilde\Delta$. Moreover, since $\widetilde\Phi(v_k), \widetilde v_k'\in\cV_2^{(k)}$ and by Lemma \ref{lemm: sigma},
\begin{align*}
\widetilde\Phi(v_k)-\widetilde v'_k\in\begin{cases}
\bQ f_k, & k \text{: odd},\\
0, & k \text{: even},
\end{cases}
\end{align*}
and we have the claim if $k$ is even. Also, if $k$ is odd, there exists $\beta'_k\in\bQ$ such that
\begin{align*}
\widetilde\Phi(v_k)=\widetilde v'_k+(-1)^{k-1}\beta'_kf_k,
\end{align*}
then, we have the claim for any positive integer $k$.
\epf

\lem \label{lem: 5}
For each positive integer $k\ge1$, we have
\begin{align}
&\widetilde\Phi\left(\frac{\pi^{\mathfrak{m}}}{2}\sum_{\mathbf{l}\in\mathbb{I}_{k-1}}\frac{(\logm2)^{l_1}}{l_1!}\prod_{j=2}^{|\mathbf{l}|}\frac{1}{l_j!}\left\{\frac{1-2^{1-j}}{j}\zetam(j)\right\}^{l_j}\right)\\
&=\left((-1)^{k-1}\sum_{(\mathbf{l};l')\in \widetilde\bI_k^{\od}}\alpha_{l'}\tau^{l'}\cdot\underset{j\ge1 \text{: odd}}{\Shu}(\beta_jf_j)^{l_j}\right)\otimes\frac{1}{\sqrt{-1}}\notag \label{eq: maininU}
\end{align}
where the coefficients $\alpha_j,\beta_j\in\bQ$ are defined by
\begin{align*}
\alpha_{n}&=\frac{1}{4}\sum_{\mathbf{l}\in\bI_{n-1}^{\ev}}\prod_{j=2\text{: even}}^{n-1}\frac{1}{l_{j}!}\cdot\left(-\frac{(1-2^{1-j})B_{j}}{2j\cdot j!}\right)^{l_{j}},\\
\beta_{n}&=\frac{2^{n}}{n}
\end{align*}
for each positive odd number $j\ge1$, and put
\begin{align*}
\begin{array}{ll}
\widetilde\bI_k^{\od}=\bigsqcup_{l'\ge1\text{: odd}}\{(\bl;l')~|~\bl\in\mathbb{I}_{k-l'}\},&\\[3mm]
\bI_0^{\ev}=\{(0,0)\},&\\[3mm]
\bI_k^{\ev}=\{(l_1,\dots,l_s)\in\bI_k~|~l_j=0~\text{if }j\text{: odd}\}, & k\ge1
\end{array}
\end{align*}
for $k\in\bZge{0}$.
\elem
\pf
Note that by Proposition \ref{prop: gk}, we have
\renewcommand{\arraystretch}{1.5}
\begin{align*}
g_j=\begin{cases}
-\frac{1}{2}\log^\mathfrak{m}2, & j=1,\\
2^{-2j+1}(1-2^{j-1})\zetam(j), & j>1.
\end{cases}
\end{align*}
\renewcommand{\arraystretch}{1}
If $j$ is odd, we have $g_j=\Phi^{-1}(f_j)$ and
\begin{align*}
\begin{array}{ll}
\widetilde\Phi(\logm2)=-2\widetilde\Phi(g_1)=-2f_1=-\beta_1f_1, & j=1,\\
\widetilde\Phi(\zetam(j))=\frac{2^{2j-1}}{1-2^{j-1}}\widetilde\Phi(g_j)=\frac{2^{2j-1}}{1-2^{j-1}}f_j=-\frac{j}{1-2^{1-j}}\beta_jf_j, & j\ge3.
\end{array}
\end{align*}
On the other hands, if $j$ is even, by Proposition \ref{prop: zetamev} we have
\begin{align*}
\widetilde\Phi(\zetam(j))=-\frac{B_j}{2j!}\tau^j.
\end{align*}
Then, by separating the production by the parity of $j$, we have
\begin{align*}
&\widetilde\Phi\left(\frac{(\pi\sqrt{-1})^{\mathfrak{m}}}{2}\sum_{\mathbf{l}\in\mathbb{I}_{k-1}}\frac{(\logm2)^{l_1}}{l_1!}\prod_{j=2}^{|\mathbf{l}|}\frac{1}{l_j!}\left\{\frac{1-2^{1-j}}{j}\zetam(j)\right\}^{l_j}
\right)\\
&=\frac{\tau}{4}\sum_{\mathbf{l}\in\mathbb{I}_{k-1}}\frac{(-\beta_1f_1)^{l_1}}{l_1!}\prod_{j\ge3\text{: odd}}\frac{1}{l_j!}\left(-\beta_jf_j\right)^{l_j}\cdot\prod_{j\ge2\text{: even}}\frac{1}{l_j!}\left(-\frac{(1-2^{1-j})B_j}{2j\cdot j!}\cdot\tau^j\right)^{l_j}\\
&=(-1)^{k-1}\frac{\tau}{4}\sum_{\mathbf{l}\in\mathbb{I}_{k-1}}\frac{(\beta_1f_1)^{l_1}}{l_1!}\prod_{j\ge3\text{: odd}}\frac{1}{l_j!}\left(\beta_jf_j\right)^{l_j}\cdot\prod_{j\ge2\text{: even}}\frac{1}{l_j!}\left(-\frac{(1-2^{1-j})B_j}{2j\cdot j!}\cdot\tau^j\right)^{l_j}\\
&=(-1)^{k-1}\frac{\tau}{4}\sum_{\mathbf{l}\in\mathbb{I}_{k-1}}\prod_{j\ge1\text{: odd}}\frac{1}{l_j!}\left(\beta_jf_j\right)^{l_j}\cdot\prod_{j\ge2\text{: even}}\frac{1}{l_j!}\left(-\frac{(1-2^{1-j})B_j}{2j\cdot j!}\right)^{l_j}\cdot\tau^{jl_j}.
\end{align*}
Now, on the right hand side of the equation above, we change the order of the sum: calculating the terms whose $l'=\sum_{j\ge2\text{: even}}jl_j+1$ is constant. Then, since $\prod_{j\ge2\text{: even}}\tau^{jl_j}=\tau^{l'-1}$, we have
\begin{align*}
&\widetilde\Phi\left(\frac{(\pi\sqrt{-1})^{\mathfrak{m}}}{2}\sum_{\bl\in\mathbb{I}_{k-1}}\frac{(\logm2)^{l_1}}{l_1!}\prod_{j=2}^{|\mathbf{l}|}\frac{1}{l_j!}\left\{\frac{1-2^{1-j}}{j}\zetam(j)\right\}^{l_j}
\right)\\
&=\frac{(-1)^{k-1}\tau}{4}\sum_{\bl\in\bI_{k-1}^{\od}}\prod_{j\ge1 \text{: odd}}\frac{1}{l_j!}(\beta_j f_j)^{l_j}\\
&~~+\frac{(-1)^{k-1}}{4}\sum_{l'\ge3 \text{: odd}}\sum_{(l_2,l_4,\dots)\in\bI_{l'-1}^{\ev}}\sum_{(l_1,l_3,\dots)\in\mathbb{I}_{k-l'}^{\od}}\prod_{j\ge1\text{: odd}}\frac{1}{l_j!}\left(\beta_jf_j\right)^{l_j}\cdot\prod_{j\ge2\text{: even}}\frac{1}{l_j!}\left(-\frac{(1-2^{1-j})B_j}{2j\cdot j!}\right)^{l_j}\cdot\tau^{l'}\\
&=\frac{(-1)^{k-1}}{4}\sum_{l'\ge1 \text{: odd}}\sum_{(l_1,l_3,\dots)\in\mathbb{I}_{k-l'}^{\od}}\prod_{j\ge1\text{: odd}}\frac{1}{l_j!}\left(\beta_jf_j\right)^{l_j}\cdot\sum_{(l_2,l_4,\dots)\in\bI_{l'-1}^{\ev}}\prod_{j\ge2\text{: even}}\frac{1}{l_j!}\left(-\frac{(1-2^{1-j})B_j}{2j\cdot j!}\right)^{l_j}\cdot\tau^{l'}\\
&=(-1)^{k-1}\sum_{l'\ge1 \text{: odd}}\sum_{(l_1,l_3,\dots)\in\mathbb{I}_{k-l'}^{\od}}\prod_{j\ge1\text{: odd}}\frac{1}{l_j!}\left(\beta_jf_j\right)^{l_j}\cdot\alpha_{l'}\tau^{l'}\\
&=(-1)^{k-1}\sum_{(\mathbf{l};l')\in \widetilde\bI_k^{\od}}\alpha_{l'}\tau^{l'}\underset{j\ge1 \text{: odd}}{\Shu}(\beta_jf_j)^{l_j}.
\end{align*}
\epf

\subsection{Proof of Theorem \ref{theo: main2}}
\pf[of Theorem \ref{theo: main2}]
By Lemma \ref{lem: 5}, it is sufficient if we show
\begin{align}
\widetilde{\Phi}(u_k)=(-1)^{k-1}\sum_{(\mathbf{l};l')\in \widetilde\bI_k^{\od}}\alpha_{l'}\tau^{l'}\underset{j\ge1 \text{: odd}}{\Shu}(\beta_jf_j)^{l_j}
\end{align}
for each $k\in\bZge{1}$. In addition, we prove at the same time that $\beta_{k-1}=\beta_{k-1}'$ for each even number $k\ge2$ ($\beta_{k-1}$ is given in Lemma \ref{lem: 5} and $\beta_{k-1}'$ is given in Lemma \ref{lemm: 4}). We prove these equations by induction on $k$. First, when $k=1$, $\widetilde\Phi(u_1)\in\cU_4^{(1)}=\bQ\tau\oplus\bQ\log^{\mathfrak{m}}2$. On the other hands, since
\begin{align*}
\per(u_k)=\frac{1}{2}(\mathrm{I}_\dch(0;\sqrt{-1};1)-\mathrm{I}_\dch(0;-\sqrt{-1};1))=\frac{\pi\sqrt{-1}}{2}
\end{align*}
and $\alpha_1=\frac{1}{4}$, we have $\widetilde\Phi(u_1)=\alpha_1\tau$, which yields the case of $k=1$ in the equation (\ref{eq: maininU}). When $k>1$, we have
\begin{align}
\widetilde\Delta(u_k)&=-\sum_{r=1}^{k-1}v'_{k-r}\otimes u_r, \notag\\
\Delta(u_k)&=-\sum_{r=1}^{k-1}v'_{k-r}\otimes u_r+1\otimes u_k+\rho(u_k)\otimes1, \label{eq: 3}
\end{align}
by the equation (\ref{eq: coaction_u}) in Lemma \ref{lemm: coaction_uv}. On the other hands, if we put
\begin{align*}
\widetilde u_{k}&:=(-1)^{k-1}\sum_{(\mathbf{l};l')\in\widetilde\bI_k^{\od}}\alpha_{l'}\tau^{l'}\underset{j\ge1 \text{: odd}}{\Shu}(\beta_jf_j)^{l_j}\in\cU_4^{(k)},\\
\widetilde v_{k}&:=(-1)^{k-1}\sum_{\mathbf{l}\in \bI_k^{\od}}\underset{j\ge1 \text{: odd}}{\Shu}(\beta_jf_j)^{l_j}\in\cV_4^{(k)},\\
\widetilde v_{m}'&:=\Phi(v_m)=(-1)^{m-1}\sum_{\mathbf{l}\in \bI_m^{\od}}\underset{j\ge1 \text{: odd}}{\Shu}(\beta_j'f_j)^{l_j}\in\cV_4^{(m)},
\end{align*}
then, we have
\begin{align*}
\Delta(\widetilde u_{k})&=(-1)^{k-1}\sum_{(\mathbf{l};l')\in \widetilde\bI_k^{\od}}\Delta(\alpha_{l'}\tau^{l'})\cdot\underset{j\ge1 \text{: odd}}{\Shu}\Delta((\beta_jf_j)^{l_j})\\
&=(-1)^{k-1}\sum_{(\mathbf{l};l')\in \widetilde\bI_k^{\od}}(1\otimes \alpha_{l'}\tau^{l'})\cdot\underset{j\ge1 \text{: odd}}{\Shu}\sum_{m_j=0}^{l_j}(\beta_jf_j)^{l_j-m_j}\otimes (\beta_jf_j)^{m_j}\\
&=(-1)^{k-1}\sum_{(\mathbf{l};l')\in \widetilde\bI_k^{\od}}\sum_{m_1+n_1=l_1}\dots\sum_{m_{|\mathbf{l}|}+n_{|\mathbf{l}|}=l_{|\mathbf{l}|}}\left(\underset{j\ge1 \text{: odd}}{\Shu}(\beta_jf_j)^{m_j}\right)\otimes \left(\alpha_{n'}\tau^{n'}\underset{j\ge1 \text{: odd}}{\Shu}(\beta_jf_j)^{n_j}\right).
\end{align*}
Now, we put $m=m_1+3m_3+\dots+|\bl|m_{|\bl|}$ and $n=n_1+3n_3+\dots+|\bl|n_{|\bl|}=k-m$, and change the order of the sum: calculating the sum of the terms for each $m=0,\dots,k$. Then, we have
\begin{align*}
\Delta(\widetilde u_{k})&=-\sum_{m+n=k}\left((-1)^{m-1}\sum_{\mathbf{m}\in \bI_m^{\od}}\underset{j\ge1 \text{: odd}}{\Shu}(\beta_jf_j)^{m_j}\right)\otimes\left((-1)^{n-1}\sum_{(\mathbf{n};n')\in\widetilde\bI_n^{\od}}\alpha_{l'}\tau^{l'}\underset{j\ge1 \text{: odd}}{\Shu}(\beta_jf_j)^{n_j}\right).
\end{align*}
Now, by the induction hypothesis, we have
\begin{align*}
\widetilde\Phi(v_m)&=(-1)^{m-1}\sum_{\mathbf{m}\in \bI_m^{\od}}\underset{j\ge1 \text{: odd}}{\Shu}(\beta_jf_j)^{m_j}\\
\widetilde\Phi(u_n)&=(-1)^{n-1}\sum_{(\mathbf{n};n')\in\widetilde\bI_n^{\od}}\alpha_{l'}\tau^{l'}\underset{j\ge1 \text{: odd}}{\Shu}(\beta_jf_j)^{n_j}
\end{align*}
for each $0<m,n<k$. Therefore,
\begin{align}
\label{eq: 4}
\Delta(\widetilde u_{k})=-\sum_{\substack{m+n=k\\m,n>0}}\widetilde\Phi(v_m)\otimes\widetilde\Phi(u_n)+1\otimes\widetilde u_k+\rho(\widetilde u_k)\otimes1.
\end{align}
Moreover, by the induction hypothesis, $\beta_m=\beta_m'$ for each $m<k-1$, and we have
\begin{align*}
\widetilde v_m'-\widetilde v_m=\begin{cases}
(-1)^{k-2}(\beta'_{k-1}-\beta_{k-1})f_{k-1}, & m=k-1 \text{: even},\\
0, & m<k-1 \text{ or } m=k-1 \text{: odd}.
\end{cases}
\end{align*}
Furthermore, by the equations (\ref{eq: 3}) and (\ref{eq: 4}), we have
\begin{align} \label{eq: 5}
\Delta(\widetilde\Phi(u_k)-\widetilde u_k)&=(\widetilde\Phi\otimes\widetilde\Phi)(\Delta(u_k))-\Delta(\widetilde u_k)\\ \notag
&=-\sum_{r=1}^{k-1}\widetilde\Phi(v'_{k-r})\otimes \widetilde\Phi(u_r)+1\otimes \widetilde\Phi(u_k)+\rho(\widetilde\Phi(u_k))\otimes1\\ \notag
&~~+\sum_{\substack{m+n=k\\m,n>0}}\widetilde\Phi(v_m)\otimes\widetilde\Phi(u_n)-1\otimes\widetilde u_k-\rho(\widetilde u_k)\otimes1\\ \notag
&=1\otimes(\widetilde\Phi(u_k)-\widetilde u_k)+(\widetilde\Phi(u_k)-\widetilde u_k)\otimes1-(\widetilde v'_{k-1}-\widetilde v_{k-1})\otimes u_1\\ \notag
&=\begin{cases}
1\otimes(\widetilde\Phi(u_k)-\widetilde u_k)+(\widetilde\Phi(u_k)-\widetilde u_k)\otimes1-(-1)^{k}(\beta'_{k-1}-\beta_{k-1})\otimes \alpha_1\tau, & k \text{: even},\\
1\otimes(\widetilde\Phi(u_k)-\widetilde u_k)+(\widetilde\Phi(u_k)-\widetilde u_k)\otimes1, & k \text{: odd}.
\end{cases}
\end{align}

Now, if $k$ is odd, $\widetilde\Phi(u_k)-\widetilde u_k\in\Ker\widetilde\Delta$. Here, by considering $\sigma$ action, $\eta$ is changed to $-\eta$, and $\eta_0$ is changed to $\eta_0$. Then, we have ${}^\sigma u_k={}^\sigma\Im(0;\eta\eta_0^{k-1};1)=-u_k$. Also, by Lemma \ref{lemm: 1}, we have ${}^\sigma\widetilde{u}_k=-\widetilde{u}_k,$ then $$\widetilde\Phi(u_k)-\widetilde u_k\in\left(\cU_4^{(k)}\right)^{\sigma,-}.$$ In particular, by Lemma \ref{lemm: sigma}, there exists $\alpha'\in\bQ$ such that $$\widetilde\Phi(u_k)-\widetilde u_k=\alpha'\tau^k.$$ Furthermore, by Proposition \ref{prop: explicit}, the equations $\per(u_k)=\per(\widetilde u_k)$ and $\per(\tau^k)\ne0$ hold, then we have $\alpha'=0$, and the claim follows.

On the other hands, if $k$ is even, we have$$\Delta((\beta'_{k-1}-\beta_{k-1})\cdot\alpha_1\tau)=(\beta'_{k-1}-\beta_{k-1})\otimes \alpha_1\tau+1\otimes(\beta'_{k-1}-\beta_{k-1})\cdot\alpha_1\tau.$$ Also, the equation (\ref{eq: 5}) yields
\begin{align*}
\Delta(u)=1\otimes u+u\otimes1,
\end{align*}
where $u:=\widetilde\Phi(u_k)-\widetilde u_k+(-1)^k(\beta'_{k-1}-\beta_{k-1})\cdot\alpha_1\tau$. Therefore, by Lemma \ref{lemm: kernelH}, we have
\begin{align*}
u=\widetilde\Phi(u_k)-\widetilde u_k+(-1)^k(\beta'_{k-1}-\beta_{k-1})\cdot\alpha_1\tau\in\bQ\tau^k\oplus\bQ f_k.
\end{align*}
Here, the equation $\beta'_{k-1}=\beta_{k-1}$ holds since $k>1$. Moreover, since ${}^\sigma u_k=-u_k, {}^\sigma\widetilde{u}_k=-\widetilde{u}_k$, there exists $\beta'\in\bQ$ such that $$\widetilde\Phi(u_k)-\widetilde u_k=\beta' f_k$$ by Lemma \ref{lemm: sigma}. Furthermore, by Proposition \ref{prop: explicit} yields the equations $\per(u_k)=\per(\widetilde u_k)$ and $\per(f_k)\ne0$. Therefore, we have $\beta'=0$, and the claim follows.

\epf

\end{document}